\documentclass[12pt]{article}

\usepackage{latexsym}
\usepackage{calc}
\usepackage{amssymb, amsmath, amsfonts, listings, mathrsfs} 
\usepackage{xcolor} 
\usepackage[normalem]{ulem}

\usepackage{graphicx}
\usepackage[labelfont=bf]{caption}

\newcounter{hours}\newcounter{minutes}

\newcommand{\stkout}[1]{\ifmmode\text{\sout{\ensuremath{#1}}}\else\sout{#1}\fi}

\def\Argmin{\mathop{\rm Argmin}}

\def\conv{{\rm Conv \,}}

\def\nr{\par \noindent}

\def\Def{\stackrel{\mathrm{def}}{=}}

\def\vf{\varphi}
\def\dom{{\rm dom \,}}
\def\beq{\begin{equation}}
\def\eeq{\end{equation}}

\def\R{\mathbb{R}}
\def\E{\mathbb{E}}
\def\S{\mathbb{S}}
\def\D{\mathscr{D}}

\def\BI{\begin{itemize}}
\def\EI{\end{itemize}}

\newcommand{\refLE}[1]{\ensuremath{\stackrel{(\ref{#1})}{\leq}}}
\newcommand{\refEQ}[1]{\ensuremath{\stackrel{(\ref{#1})}{=}}}
\newcommand{\refGE}[1]{\ensuremath{\stackrel{(\ref{#1})}{\geq}}}

\newtheorem{theorem}{Theorem}
\newtheorem{lemma}{Lemma}
\newtheorem{corollary}{Corollary}

\newtheorem{proposition}{Proposition}
\newtheorem{assumption}{Assumption}
\newtheorem{definition}{Definition}

\newtheorem{example}{Example}
\newtheorem{remark}{Remark}
\newcommand{\proof}{\bf Proof: \rm \nr}
\newcommand{\qed}{\hfill $\Box$ \nr \medskip}
\newcommand{\half}{\mbox{${1 \over 2}$}}

\def\ba{\begin{array}}
\def\ea{\end{array}}
\def\beann{\begin{eqnarray*}}
\def\eeann{\end{eqnarray*}}
\def\bea{\begin{eqnarray}}
\def\eea{\end{eqnarray}}

\def\BT{\begin{theorem}}
\def\ET{\end{theorem}}
\def\BL{\begin{lemma}}
\def\EL{\end{lemma}}
\def\BC{\begin{corollary}}
\def\EC{\end{corollary}}
\def\BE{\begin{example}}
\def\EE{\end{example}}
\def\BD{\begin{definition}}
\def\ED{\end{definition}}
\def\BR{\begin{remark}}
\def\ER{\end{remark}}
\def\BAS{\begin{assumption}}
\def\EAS{\end{assumption}}
\def\BI{\begin{itemize}}
\def\EI{\end{itemize}}
\def\BP{\begin{proposition}}
\def\EP{\end{proposition}}

\def\BMP{\begin{minipage}{9.5cm}}
\def\EMP{\end{minipage}}
\def\MPT{\begin{minipage}{11.5cm}}
\def\EPT{\end{minipage}}

\def\la{\langle}
\def\ra{\rangle}

\def\QF{\hspace{5ex} \Box}
\def\QR{\hfill \Box}

\title{
\vspace{10mm} 
\textbf{Affine-invariant
contracting-point methods \\ for Convex Optimization}
\thanks{Research results presented in this paper were obtained
in the framework of ERC Advanced Grant~788368.} }
\author{Nikita Doikov
\thanks{Institute of Information and Communication Technologies,
Electronics and Applied Math. (ICTEAM), Catholic
University of Louvain (UCL). E-mail:
{Nikita.Doikov@uclouvain.be}. ORCID: 0000-0003-1141-1625.} \and Yurii Nesterov
\thanks{Center for Operations Research and Econometrics (CORE),
Catholic University of Louvain (UCL). \newline E-mail:
{Yurii.Neterov@uclouvain.be}. ORCID: 0000-0002-0542-8757.}
}

\date{\today}

\begin{document}
\maketitle

\begin{abstract}
In this paper, we develop new affine-invariant algorithms
for solving composite convex minimization problems
with bounded domain. We present a general
framework of Contracting-Point methods, which solve at each iteration an auxiliary subproblem restricting the smooth part of the objective function
onto contraction of the initial domain.
This framework provides us with a systematic way 
for developing optimization methods of different order,
endowed with the global complexity bounds.
We show that using an appropriate affine-invariant smoothness condition,
it is possible to implement one iteration of the Contracting-Point method
by one step of the pure tensor method of degree $p \geq 1$.
The resulting global rate of convergence in functional residual is then
${\cal O}(1 / k^p)$, where $k$ is the iteration counter. 
It is important that all constants in our bounds are {\em affine-invariant}.
For $p = 1$, our scheme recovers well-known Frank-Wolfe algorithm,
providing it with a new interpretation by a general perspective 
of tensor methods. Finally, within our framework, 
we present efficient implementation and total complexity analysis
of the inexact second-order scheme $(p = 2)$,
called Contracting Newton method. It can be seen as a proper implementation of the {\em trust-region idea}.
Preliminary numerical results confirm its good practical performance both in the number of iterations, and in computational time.

\end{abstract}

\vspace{3ex}\noindent

\thispagestyle{empty}

\newpage\setcounter{page}{1}

\section{Introduction}
\setcounter{equation}{0}

\vspace{1ex}\noindent
{\bf Motivation.} In the last years, we can see an
increasing interest to new frameworks for derivation
and justification different methods for Convex Optimization,
provided with a worst-case complexity analysis (see, for
example, \cite{bauschke2016descent, lu2018relatively, cartis2018global, nesterov2018complexity, nesterov2019implementable, gasnikov2019near,
doikov2019contracting, kamzolov2020optimal, nesterov2020superfast,
nesterov2020inexact}). It appears that the accelerated proximal
tensor methods \cite{baes2009estimate, nesterov2019implementable} can be naturally
explained through the framework of high-order
proximal-point schemes \cite{nesterov2020inexact} requiring solution
of nontrivial auxiliary problem at every iteration.

This possibility serves as a departure point for the
results presented in this paper. Indeed, the main drawback
of proximal tensor methods consists in necessity of using
a fixed Euclidean structure for measuring distances
between points. However, the multi-dimensional Taylor
polynomials are defined by directional derivatives, which
are affine-invariant objects. Can we construct a family of
tensor methods, which do not depend on the choice of the
coordinate system in the space of variables? The results
of this paper give a positive answer on this question.

Our framework extends the initial results presented in
\cite{nesterov2018complexity} and in \cite{doikov2020convex}. In
\cite{nesterov2018complexity}, it was shown that the classical Frank-Wolfe
algorithm can be generalized onto the case of the
composite objective function \cite{nesterov2013gradient} using a
contraction of the feasible set towards the current test
point. This operation was used there also for justifying a
second-order method with contraction, which looks similar
to the classical trust-region methods \cite{conn2000trust}, but with
asymmetric trust region.
The convergence rates for the second-order methods with contractions were significantly improved in~\cite{doikov2020convex}.
In this paper, we extend the
contraction technique onto the whole family of tensor
methods. However, in the
vein of \cite{nesterov2020inexact}, we start first from analyzing a
conceptual scheme solving at each iteration an auxiliary
optimization problem formulated in terms of the initial
objective function.

The results of this work can be also seen 
as an affine-invariant counterpart of
Contracting Proximal Methods from~\cite{doikov2019contracting}.
In the latter algorithms, one need to fix the prox function which is 
suitable for the geometry of the problem, in advance.
The parameters of the problem class are also usually required.
The last but not least, all methods from this work are universal
and parameter-free.

\vspace{1ex}\noindent
{\bf Contents.} 
The paper is organized as follows.

In Section~\ref{sc-CPM}, we present 
a general framework of Contracting-Point methods.
We provide two conceptual variants of our scheme
for different conditions of inexactness 
for the solution of the subproblem: 
using a point with small residual in the function value,
and using a stronger condition which involves the gradients.
For both schemes we establish global bounds for the functional residual
of the initial problem.
These bounds lead to global convergence guarantees under a suitable choice of
the parameters.
For the scheme with the second condition of inexactness,
we also provide a computable accuracy certificate.
It can be used to estimate the functional residual directly within the method.

Section~\ref{sc-Aff} contains smoothness conditions, which are useful to analyse
affine-invariant high-order schemes. We present some basic inequalities and examples,
related to the new definitions.

In Section~\ref{sc-Ten}, we show how to implement one iteration of our methods 
by computing an (inexact) affine-invariant tensor step.
For the methods of degree $p \geq 1$, we establish global convergence 
in the functional residual
of the order ${\cal O}(1 / k^p)$, where $k$ is the iteration counter.
For $p = 1$, this recovers well-known result about global convergence of 
the classical Frank-Wolfe algorithm~\cite{frank1956algorithm,nesterov2018complexity}.
For $p = 2$, we obtain Contracting-Domain Newton Method from~\cite{doikov2020convex}.
Thus, our analysis also extends the results from these works to the case, 
when the corresponding subproblem is solved inexactly.

In Section~\ref{sc-Newton}, we present two-level optimization scheme,
called Inexact Contracting Newton Method.
This is implementation 
of the inexact second-order method, via computing its steps by the first-order
Conditional Gradient Method.
For the resulting algorithm, 
we establish global complexity ${\cal O}(1 /  \varepsilon^{1/2})$
calls of the \textit{smooth part oracle} 
(computing gradient and Hessian of the smooth part of the objective),
and ${\cal O}(1 / \varepsilon)$ calls
of the \textit{linear minimization oracle} of the composite part,
where $\varepsilon > 0$ is the required accuracy in the functional residual.
Additionally, we address effective implementation of our method for optimization
over the standard simplex.

Section~\ref{sc-Numerical} contains numerical experiments.

In Section~\ref{sc-Discussion}, we discuss our results and highlight some
open questions for the future research.

\vspace{1ex}\noindent
{\bf Notation.} In what follows we denote by $\E$ a
finite-dimensional real vector space, and by $\E^*$ its
dual space, which is a space of linear functions on $\E$.
The value of function $s \in \E^{*}$ at point $x \in \E$
is denoted by $\la s, x \ra$.

For a smooth function $f: \dom f \to \R$, where $\dom f
\subseteq \E$, we denote by $\nabla f(x)$ its gradient and
by $\nabla^2 f(x)$ its Hessian, evaluated at point $x \in
\dom f \subseteq \E$. Note that
$$
\ba{rcl}
\nabla f(x) & \in & \E^{*}, \qquad \nabla^2 f(x) h \;\;
\in \;\; \E^{*},
\ea
$$
for all $x \in \dom f$ and $h \in \E$. For $p \geq 1$, we
denote by $ D^p f(x)[h_1, \dots, h_p] $ the $p$th
directional derivative of $f$ along directions $h_1,
\dots, h_p \in \E$. Note that $D^p f(x)$ is a $p$-linear
symmetric form on $\E$. If $h_i = h$ for all $1 \leq i
\leq p$, a shorter notation $D^p f(x)[h]^p$ is used. For
its gradient in $h$, we use the following notation:
$$
\ba{rcl}
 D^p f(x)[h]^{p - 1} \; \Def \; {1 \over p}\nabla D^p f(x)[h]^p  & \in & \E^{*}, \qquad h \in \E.
\ea
$$
In particular, $D^1 f(x)[h]^{0} \equiv \nabla f(x)$,
and $D^2 f(x)[h]^{1} \equiv \nabla^2 f(x)h$.

\section{Contracting-point methods}
\label{sc-CPM}

Consider the following composite minimization problem
\beq \label{MainProblem}
\ba{rcl}
F^* \; \Def \; \min\limits_{x \in \dom \psi} \Big[F(x) & =
& f(x) + \psi(x) \Big],
\ea
\eeq
where $\psi : \E \to \R \cup \{ +\infty \}$ is a {\em
simple} proper closed convex function with bounded domain,
and function $f(x)$ is convex and $p$ ($\geq 1$) times
continuously differentiable at every point $x \in \dom
\psi$.

In this section, we propose a conceptual optimization
scheme for solving~\eqref{MainProblem}. At each step of
our method, we choose a contracting coefficient $\gamma_k
\in (0, 1]$ restricting the nontrivial part of our
objective $f(\cdot)$ onto a \textit{contracted domain}. At
the same time, the domain for the composite part remains
unchanged.

Namely, at point $x_k \in \dom \psi$, define
$$
\ba{rcl}
S_k(y) & \Def & \gamma_k \psi\bigl( x_k +
\frac{1}{\gamma_k}(y - x_k) \bigr),
\quad y = x_k + \gamma_k(v - x_k),
\quad v \in \dom \psi.
\ea
$$
Note that $S_k(y) = \gamma_k \psi(v)$. Consider the
following {\em exact} iteration:
\beq \label{ExactMethod}
\boxed{
\ba{rcl}
v^*_{k + 1} & \in & \Argmin\limits_{v} \Bigl\{ f(y) +
S_k(y) : \; y = (1 - \gamma_k )x_k +
\gamma_k v, \, v \in \dom \psi \Bigr\},\\
x^*_{k+1} & = & (1 - \gamma_k )x_k + \gamma_k v^*_{k+1}.
\ea
}
\eeq
Of course, when $\gamma_k = 1$, exact step from~\eqref{ExactMethod}
solves the initial problem.
However, we are going to look at the \textit{inexact} minimizer.
In this case, the choice of $\{ \gamma_k \}_{k \geq 0}$ should take into account
the efficiency of solving the auxiliary subproblem.

Denote by $F_k(\cdot)$ the objective in the auxiliary
problem \eqref{ExactMethod}, that is
$$
\ba{rcl}
F_k(y) & \Def & f(y) + S_k(y), \quad y \; = \;
(1-\gamma_k) x_k + \gamma_k v, \quad v \in \dom \psi.
\ea
$$
We are going to use the point $\bar x_{k + 1} =
(1-\gamma_k) x_k + \gamma_k \bar v_{k+1}$ with $\bar
v_{k+1} \in \dom \psi$ having a \textit{small residual} in
the function value:
\beq \label{Inex1}
\ba{rcl}
F_k(\bar x_{k + 1}) - F_k(x_{k+1}^*) & \leq & \delta_{k +
1},
\ea
\eeq
with some fixed $\delta_{k + 1} \geq 0$.
\BL\label{lm-Iter}
For all $k \geq 0$ and $v \in \dom \psi$, we have
\beq\label{eq-Iter}
\ba{rcl}
F(\bar x_{k+1}) & \leq & (1-\gamma_k) F(x_k) + \gamma_k
F(v) + \delta_{k + 1}.
\ea
\eeq
\EL
\proof
Indeed, for any $v \in \dom \psi$, we have
$$
\ba{rcl}
F_k(\bar x_{k+1}) & \refLE{Inex1} & F_k(x_{k+1}^*) +
\delta_{k+1}\\
\\
& \refLE{ExactMethod} & f((1 - \gamma_k )x_k + \gamma_k v)
+ S_k((1 - \gamma_k )x_k + \gamma_k v) +
\delta_{k+1}\\
\\
& \leq & (1-\gamma_k) f(x_k) + \gamma_k f(v) + \gamma_k
\psi(v) + \delta_{k+1}.
\ea
$$
Therefore,
$$
\ba{rcl}
F(\bar x_{k+1}) & = & F_k(\bar x_{k+1}) + \psi(\bar
x_{k+1}) - \gamma_k\psi(\bar v_{k + 1}) \\
\\
& \leq & (1-\gamma_k) f(x_k) + \gamma_k F(v) +
\delta_{k+1} + \psi(\bar x_{k+1}) - \gamma_k \psi(\bar v_{k+1})\\
\\
& \leq & (1-\gamma_k) F(x_k) + \gamma_k F(v) +
\delta_{k+1}. \QR
\ea
$$

Let us write down our method in an algorithmic form.
\beq\label{met-CPoint1}
\ba{|c|}
\hline \\
\mbox{\bf Conceptual Contracting-Point Method, I}\\
\\
\hline \\
\ba{l}
\mbox{{\bf Initialization.} Choose $x_0 \in \dom \psi$.}\\
\\
\mbox{\bf Iteration $k \geq 0$.}\\
\\
\mbox{1: Choose $\gamma_k \in (0,1]$.}\\
\\
\mbox{2: For some $\delta_{k + 1} \geq 0$, find
$\bar{x}_{k+1}$ satisfying (\ref{Inex1}).
    }\\
\\
\mbox{3: If $F(\bar x_{k + 1}) \leq F(x_k)$, then set $x_{k + 1} = \bar x_{k + 1}$.
Else choose $x_{k+1} = x_k$.}\\
\ea\\
\\
\hline
\ea
\eeq

In Step 3 of this method, we add a simple test for
ensuring monotonicity in the function value. This step is
optional.

It is more convenient to describe the rate of convergence
of this scheme with respect to another sequence of
parameters. Let us introduce an arbitrary sequence of
positive numbers $\{ a_k  \}_{k \geq 1}$ and denote $A_k
\Def \sum_{i = 1}^k a_i$. Then, we can define the
contracting coefficients as follows
\beq\label{def-Gamma}
\ba{rcl}
\gamma_k & \Def & \frac{a_{k + 1}}{A_{k + 1}}.
\ea
\eeq
\BT
For all points of sequence $\{ x_k \}_{k \geq 0}$,
generated by process~\eqref{met-CPoint1}, we have the
following relation:
\beq \label{CPoint1Conv}
\ba{rcl}
A_kF(x_k) & \leq & A_k F^{*} + B_k,
\quad \text{with} \quad
B_k \; \Def \; \sum_{i = 1}^k A_i \delta_i.
\ea
\eeq
\ET
\proof
Indeed, for $k = 0$, we have $A_k = 0$, $B_k = 0$. Hence,
(\ref{CPoint1Conv}) is valid. Assume it is valid for some $k
\geq 0$. Then
$$
\ba{rcl}
A_{k+1} F(x_{k+1}) & \stackrel{\mbox{Step 3}}{\leq} &
A_{k+1} F(\bar x_{k+1}) \; \leq \; A_{k+1} \Big(
(1-\gamma_k) F(x_k) + \gamma_k F^* + \delta_{k+1} \Big)\\
\\
& \refEQ{def-Gamma} & A_k F(x_k) + a_{k+1} F^* + A_{k+1}
\delta_{k+1} \; \stackrel{(\ref{CPoint1Conv})}{\leq} \;
A_{k+1} F^* + B_{k+1}. \QF
\ea
$$

From bound~\eqref{CPoint1Conv}, we can see, that
\beq \label{FConverg1}
\ba{rcl}
F(x_k) - F^{*} & \leq & \frac{1}{A_k}\, \sum_{i = 1}^k A_i
\delta_i, \qquad k \geq 1.
\ea
\eeq
Hence, the actual rate of convergence of method
(\ref{met-CPoint1}) depends on the growth of coefficients
$\{ A_k \}_{k \geq 1}$ {\em relatively} to the level of
inaccuracies $\{ \delta_k \}_{k \geq 1}$. Potentially,
this rate can be arbitrarily high. Since we did not assume
anything yet about our objective function, this means that
we just retransmitted the complexity of solving the
problem (\ref{MainProblem}) onto a lower level, the level
of computing the point $\bar x_{k+1}$, satisfying the
condition (\ref{Inex1}). We are going to discuss different
possibilities for that in Sections~\ref{sc-Ten} and~\ref{sc-Newton}.

Now, let us endow the method (\ref{met-CPoint1}) with a
computable \textit{accuracy certificate}. For this
purpose, for a sequence of given test points $\{ \bar x_k
\}_{k \geq 1} \subset \dom \psi$, we introduce the
following \textit{Estimating Function}
(see~\cite{nesterov2018lectures}):
$$
\ba{rcl}
\vf_k(v) & \Def & \sum\limits_{i = 1}^k a_i \bigl[  f(\bar
x_i) + \la \nabla f(\bar x_i), v - \bar x_i \ra + \psi(v)
\bigr].
\ea
$$
By convexity of $f(\cdot)$, we have $A_k F(v) \geq
\vf_k(v)$ for all $v \in \dom \psi$. Hence, for all $k
\geq 1$, we can get the following bound for the functional
residual:
\beq \label{LkDef}
\ba{rcl}
F(x_k) - F^{*} & \leq & \ell_k \; \Def \; F(x_k) -
\frac{1}{A_k}\, \vf_k^{*}, \quad \vf_k^* \; \Def \;
\min\limits_{v \in \dom \psi} \; \vf_k(v).
\ea
\eeq
The complexity of computing the value of $\ell_k$ usually
does not exceed the complexity of computing the next
iterate of our method since it requires just one call of
the \textit{linear minimization oracle}. Let us show that
an appropriate rate of decrease of the estimates $\ell_k$
can be guaranteed by sufficiently accurate steps of the
method~\eqref{ExactMethod}.

For that, we need a stronger condition on point $\bar x_{k
+ 1}$, that is
\beq \label{Inex2}
\ba{rcl}
\la \nabla f(\bar x_{k + 1}), v - \bar v_{k + 1} \ra + \psi(v)
& \geq & \psi(\bar v_{k + 1}) - \frac{1}{\gamma_k} \delta_{k + 1}, \quad v \in \dom \psi, \\
\\
\bar x_{k + 1} & = & (1-\gamma_k) x_k + \gamma_k \bar v_{k+1},
\ea
\eeq
with some $\delta_{k + 1} \geq 0$. Note that, for
$\delta_{k + 1} = 0$, condition \eqref{Inex2} ensures the
exactness of the corresponding step of
method~\eqref{ExactMethod}.

Let us consider now the following algorithm.
\beq\label{met-CPoint2}
\ba{|c|}
\hline \\
\mbox{\bf Conceptual Contracting-Point Method, II}\\
\\
\hline \\
\ba{l}
\mbox{{\bf Initialization.} Choose $x_0 \in \dom \psi$.}\\
\\
\mbox{\bf Iteration $k \geq 0$.}\\
\\
\mbox{1: Choose $\gamma_k \in (0,1]$.}\\
\\
\mbox{2: For some $\delta_{k + 1} \geq 0$, find
$\bar{x}_{k+1}$ satisfying (\ref{Inex2}).
    }\\
\\
\mbox{3: If $F(\bar x_{k + 1}) \leq F(x_k)$, then set
$x_{k + 1} = \bar x_{k + 1}$.
Else choose $x_{k+1} = x_k$.}\\
\ea\\
\\
\hline
\ea
\eeq

This scheme differs from the previous
method~\eqref{met-CPoint1} only in the characteristic
condition \eqref{Inex2} for the next test point.
\BT
For all points of the sequence $\{ x_k \}_{k \geq 0}$,
generated by the process~\eqref{met-CPoint2}, we have
\beq \label{Prox2Conv}
\ba{rcl}
\vf_k^{*} & \geq & A_k F(x_k) - B_k, \quad k \geq 0.
\ea
\eeq
\ET
\proof
For $k = 0$, relation~\eqref{Prox2Conv} is valid since
both sides are zeros. Assume that~\eqref{Prox2Conv}
holds for some $k \geq 0$. Then, for any $v \in \dom
\psi$, we have
$$
\ba{rcl}
\vf_{k + 1}(v) & \equiv & \vf_k(v) + a_{k + 1} \bigl[
f(\bar x_{k + 1})
+ \la \nabla f( \bar x_{k + 1} ), v - \bar x_{k + 1} \ra + \psi(v) \bigr] \\
\\
& \refGE{Prox2Conv} &
A_k F(x_k) - B_k + a_{k + 1} \bigl[  f(\bar x_{k + 1})
+ \la \nabla f( \bar x_{k + 1} ), v - \bar x_{k + 1} \ra + \psi(v) \bigr] \\
\\
& \overset{(*)}{\geq} & A_{k + 1} \bigl[ f(\bar x_{k + 1})
+ \la \nabla f(\bar x_{k + 1}), \frac{a_{k + 1} v + A_k
x_k}{A_{k + 1}} - \bar x_{k + 1} \ra \bigr]
+ A_k \psi(x_k) + a_{k + 1} \psi(v) - B_k \\
\\
& = & A_{k + 1} f(\bar x_{k + 1}) + a_{k + 1} \bigl[ \la
\nabla f(\bar x_{k + 1}), v - \bar v_{k + 1} \ra + \psi(v)
\bigr]
+ A_k \psi(x_k) - B_k \\
\\
& \refGE{Inex2} &
A_{k + 1} f(\bar x_{k + 1}) + a_{k + 1} \psi(\bar v_{k + 1}) + A_k \psi(x_k) - B_{k + 1} \\
\\
& \overset{(**)}{\geq} &
A_{k + 1} F(\bar x_{k + 1}) - B_{k + 1}
\;\; \stackrel{\mbox{Step 3}}{\geq} \;\; A_{k + 1} F(x_{k + 1}) - B_{k + 1}.
\ea
$$
Here, the inequalities $(*)$ and $(**)$ are justified by
convexity of $f(\cdot)$ and $\psi(\cdot)$,
correspondingly. Thus,~\eqref{Prox2Conv} is proved for all
$k \geq 0$.
\qed

Combining now \eqref{LkDef} with~\eqref{Prox2Conv}, we
obtain
\beq \label{FConverg2}
\ba{rcl}
F(x_k) - F^{*} & \leq & \ell_k \;\; \leq \;\;
\frac{1}{A_k} \, \sum_{i = 1}^k A_i \delta_i, \quad k \geq
1.
\ea
\eeq
We see that the right hand side in~\eqref{FConverg2}
is the same, as that one in~\eqref{FConverg1}.
However, this convergence is stronger, 
since it provides a bound for the accuracy certificate $\ell_k$.

\section{Affine-invariant high-order smoothness
	conditions}\label{sc-Aff}

We are going to describe efficiency of solving the
auxiliary problem in~\eqref{ExactMethod} by some
\textit{affine-invariant} characteristics of variation of
function $f(\cdot)$ over the compact convex sets. For a
convex set $Q$, define
\beq \label{def-Delta}
\ba{rcl}
\Delta^{(p)}_Q(f) & \Def & \sup\limits_{{ x, v \in Q,} \atop {t
		\in (0, 1] } }
\frac{1}{t^{p + 1}}\Big| f( x + t(v - x)) - f(x) -
\sum\limits_{i = 1}^p \frac{t^i}{i!}D^i f(x)[v - x]^i
\Big|.
\ea
\eeq
Note, that for $p = 1$ this characteristic was considered in~\cite{jaggi2013revisiting}
for the analysis of the classical Frank-Wolfe algorithm.

In many situations, it is more convenient to use an upper
bound for $\Delta^{(p)}_Q(f)$, which is a full variation of
its $(p+1)$th derivative over the set $Q$:
\beq \label{def-Var}
\ba{rcl}
{\cal V}^{(p+1)}_Q(f) & \Def & \sup\limits_{x, y, v\in Q}
\Big| 
D^{p+1}f(y)[v - x]^{p+1}
\Big|.
\ea
\eeq
Indeed, by Taylor formula, we have
$$
\ba{rl}
& \frac{1}{t^{p + 1}}\Big[f(x + t(v - x)) - f(x) -
\sum\limits_{i = 1}^p \frac{t^i}{i!}D^i f(x)[v - x]^i
\Big]\\
\\
= & {1 \over p!}\int\limits_0^1 (1-\tau)^p D^{p+1}f(x +
\tau t(v-x))[v-x]^{p+1} d \tau.
\ea
$$
Hence,
\beq\label{eq-DV}
\ba{rcl}
\Delta^{(p)}_Q(f) & \leq & {1 \over (p+1)!} {\cal V}^{(p+1)}_Q(f).
\ea
\eeq

Sometimes, in order to exploit a \textit{primal-dual}
structure of the problem, we need to work with the dual
objects (gradients), as in method~\eqref{met-CPoint2}. In
this case, we need a characteristic of variation of the
gradient $\nabla f(\cdot)$ over the set $Q$:
\beq \label{def-GammaV}
\ba{rcl}
\Gamma^{(p)}_Q(f) & \Def & \sup\limits_{{ x, y, v \in Q,} \atop
	{t
		\in (0, 1] } }
\frac{1}{t^p}\Big| \la \nabla f(x + t(v - x)) - \nabla
f(x) - \sum\limits_{i = 2}^{p - 1} \frac{t^{i - 1}}{(i -
	1)!}D^i f(x)[v - x]^{i - 1}, v - y \ra \Big|.
\ea
\eeq
Since
$$
\ba{rl}
& \frac{1}{t} \Bigl[ f(x + t(v - x)) - f(x) -
\sum\limits_{i = 1}^p \frac{t^i}{i!}D^i f(x)[v - x]^i
\Bigr] \\
\\
= & \frac{1}{t}\Bigl[ \int\limits_0^1 \la \nabla f(x +
\tau t(v - x)), t (v - x) \ra d\tau
- \sum\limits_{i = 1}^p \frac{t^i}{i!} D^i f( x )[v - x]^i  \Bigr] \\
\\
= & \int\limits_0^1 \la \nabla f( x + \tau t (v - x) ) -
\sum\limits_{i = 1}^p \frac{(\tau t)^{i - 1}}{(i - 1)!}
D^{i} f(x)[v - x]^{i - 1}, v - x\ra d\tau,
\ea
$$
we conclude that
\beq \label{eq-DG}
\ba{rcl}
\Delta^{(p)}_Q(f) & \leq & \frac{1}{p + 1}\Gamma^{(p)}_Q (f).
\ea
\eeq
At the same time, by Taylor formula, we get
\beq \label{GBound}
\ba{cl}
& \frac{1}{t^p}\Bigl[ \nabla f(x + t(v - x)) - \nabla f(x)
- \sum\limits_{i = 2}^{p - 1} \frac{t^{i - 1}}{(i - 1)!}
D^{i} f(x)[v - x]^{i - 1}
\Bigr] \\
\\
= & \frac{1}{(p - 1)!} \int\limits_0^1 (1 - \tau)^{p - 1}
D^{p + 1} f(x + \tau t(v - x) )[v - x]^{p} d \tau.
\ea
\eeq
Therefore, again we have an upper bound in terms of the
variation of $(p + 1)$th derivative, that is
$$
\ba{rcl}
\Gamma^{(p)}_Q(f) & \refLE{GBound} & \frac{1}{p!} \sup\limits_{x, y, z, v \in Q}
\Big| \la D^{p+1}f(z)[v-x]^{p}, v - y\ra
\Big| \;\; \leq \;\; \frac{2(p + 1)^p}{(p!)^2} {\cal
	V}^{(p+1)}_Q(f).
\ea
$$
See Proposition~\ref{propMultilinear} in Appendix
for the proof of the last inequality.
Hence, the value of ${\cal V}_Q^{(p + 1)}(f)$ is the biggest
one. However, in many cases it is more convenient.

\BE \label{ExampleBasic}
For a fixed self-adjoint positive-definite linear operator
$B: \E \to \E^{*}$, define the corresponding Euclidean
norm as $\|x\| := \la Bx, x \ra^{1/2}, \; x \in \E$.
Let $p$th derivative of function $f(\cdot)$ be Lipschitz
continuous with respect to this norm:
$$
\ba{rcl}
\| D^p f(x) - D^p f(y) \| & \Def & \max\limits_{h \in \E:
    \|h\| \leq 1} | (D^p f(x) - D^p f(y))[h]^p | \; \leq \;
L_p \|y - x\|,
\ea
$$
for all $x, y \in Q$.
Let $Q$ be a compact set
with diameter
$$
\ba{rcl}
\D & = & \D_{\| \cdot \|}(Q) \; \Def \; \max\limits_{x, y \in Q}
\|x - y\| \; < \; +\infty.
\ea
$$
Then, we have
$$
\ba{rcl}
{\cal V}^{(p+1)}_Q(f) & \leq & L_p \D^{p + 1}.
\ea
$$
\EE

In some situations we can obtain much better estimates.

\BE
Let $A \succeq 0$, and $f(x) = \half \la A x, x \ra$ with
$x \in \S_n \Def \{ x \in \R^n_+: \sum\limits_{i=1}^n
x^{(i)} = 1 \}$. For measuring distances in the standard
simplex, we choose $\ell_1$-norm:
$$
\ba{rcl}
\| h \| & = & \sum\limits_{i=1}^n |h^{(i)}|, \quad h \in
\R^n.
\ea
$$
In this case, $\D = \D_{|| \cdot \|}(\S_n) = 2$, and $L_1 =
\max\limits_{1 \leq i \leq n} A^{(i,i)}$. On the other
hand,
$$
\ba{rcl}
{\cal V}^{(2)}_{\S_n}(f) & = & \max\limits_{1 \leq i, j \leq
n} \la A(e_i-e_j), e_i - e_j \ra \; \leq \; \max\limits_{1
\leq i, j \leq n} [2 \la Ae_i,e_i \ra + 2 \la A e_j, e_j
\ra] \; = \; 4L_1,
\ea
$$
where $e_k$ denotes the $k$th coordinate vector in $\R^n$.
Thus, ${\cal V}^{(2)}_{\S_n} \leq L_1 \D^2$.

However, for some matrices, the value ${\cal
V}^{(2)}_{\S_n}(f)$ can be much smaller than $L_1\D^2$.
Indeed, let $A = a a^T$ for some $a \in \R^n$. Then $L_1 =
\max\limits_{1 \leq i \leq n} (a^{(i)})^2$, and
$$
\ba{rcl}
{\cal V}^{(2)}_{\S_n}(f) & = & \left[ \max\limits_{1 \leq i
\leq n} a^{(i)} - \min\limits_{1 \leq i \leq n} a^{(i)}
\right]^2,
\ea
$$
which can be much smaller than $4L_1$.
\EE

\BE
Let given vectors $a_1, \dots, a_m$ span the whole $\R^n$. Consider the objective
$$
\ba{rcl}
f(x) & = & \ln \biggl( \sum\limits_{k = 1}^m e^{\la a_k, x \ra}  \biggr),
\qquad x \in \S_n.
\ea
$$
Then, it holds (see Example~1 in~\cite{doikov2019minimizing} for the first inequality):
$$
\ba{rcl}
\la \nabla^2 f(x)h, h \ra & \leq & \max\limits_{1  \leq k, l \leq m} \la a_k - a_l, h \ra^2
\;\; \leq \;\; \max\limits_{1  \leq k, l \leq m} \|a_k - a_l\|_{\infty}^2 \|h\|_1^2, \qquad h \in \R^n.
\ea
$$
Therefore, in $\ell_1$-norm we have 
$L_1 = \max\limits_{1  \leq k, l \leq m}  \max\limits_{1 \leq i \leq n} \left[  a_k^{(i)} - a_l^{(i)} \right]^2$.
At the same time,
$$
\ba{rcl}
{\cal V}_{\S_n}^{(2)}(f)
& = & \sup\limits_{x \in \S_n} 
\max\limits_{1 \leq i, j \leq n} \la \nabla^2 f(x) (e_i - e_j), e_i - e_j \ra \\
\\
& \leq &
\max\limits_{1 \leq k, l \leq m}
\max\limits_{1 \leq i, j \leq n} 
\left[
\bigl( a_k^{(i)} - a_k^{(j)} \bigr) - 
\bigl( a_l^{(i)} - a_l^{(j)} \bigr)
\right]^2.
\ea
$$
The last expression is the maximal difference between variations of the coordinates.
It can be much smaller than $L_1 \D^2 = 4 L_1$.

Moreover, we have (see Example~1 in~\cite{doikov2019minimizing}):
$$
\ba{rcl}
|D^3 f(x)[h]^3| & \leq &
\max\limits_{1 \leq k, l \leq m} |\la a_k - a_{l}, h \ra|^3, \qquad h \in \R^n.
\ea
$$
Hence, we obtain
$$
\ba{rcl}
{\cal V}_{\S_n}^{(3)}(f) & \leq & 
\max\limits_{1 \leq k, l \leq m} \max\limits_{1 \leq i, j \leq n}
\Bigl| 
\bigl( a_k^{(i)} - a_k^{(j)} \bigr) - 
\bigl( a_l^{(i)} - a_l^{(j)} \bigr)
\Bigr|^3.
\ea
$$
\EE

\section{Contracting-point tensor methods}
\label{sc-Ten}

In this section, we show how to implement
Contracting-point methods, by using
affine-invariant tensor steps. At each iteration
of~\eqref{ExactMethod}, we approximate $f(\cdot)$ by
Taylor's polynomial of degree $p \geq 1$ around the
current point $x_k$:
$$
\ba{rcl}
f(y) & \approx & \Omega_p(f, x_k; y)
\;\; \Def \;\;
f(x_k) + \sum\limits_{i = 1}^p \frac{1}{i!} D^i f(x_k) [y - x_k]^i.
\ea
$$
Thus, we need to solve the following auxiliary problem:
\beq \label{TensorSub}
\ba{rcl}
\min\limits_{v \in \dom \psi} \Bigl\{
    M_k(y) & \Def &  \Omega_p(f, x_k; y) + S_k(y) : \;
    y = (1 - \gamma_k )x_k + \gamma_k v
\Bigr\}.
\ea
\eeq
Note that this global minimum $M_k^{*}$ is well defined since $\dom \psi$
is bounded. Let us take
$$
\ba{rcl}
\bar x_{k + 1} & = & (1 - \gamma_k
)x_k + \gamma_k \bar v_{k+1},
\ea
$$
where $\bar v_{k+1}$ is an inexact solution to
\eqref{TensorSub} in the following sense:
\beq \label{MkMin}
\ba{rcl}
M_k(\bar x_{k + 1}) - M_k^{*} & \leq & \xi_{k + 1}.
\ea
\eeq
Then, this point serves as a good candidate for the
inexact step of our method.
\BT \label{prop-TensorStep1}
Let $\xi_{k + 1} \leq c\gamma_k^{p + 1}$, for some constant $c \geq 0$. Then
$$
\ba{rcl}
F_k(\bar x_{k + 1}) - F_k^{*} & \leq & \delta_{k + 1},
\ea
$$
for $\delta_{k + 1} = (c + 2\Delta^{(p)}_{\dom \psi}(f))\gamma_k^{p + 1}$.
\ET
\proof
Indeed, for $y = x_k + \gamma_k (v - x_k)$ with arbitrary $v \in \dom \psi$, we have
$$
\ba{rcl}
F_k(y) & = & f(y) + S_k(y) \\
\\
& \refGE{def-Delta} &
\Omega_p(f, x_k; y) + S_k(y) - \Delta^{(p)}_{\dom \psi}(f)\gamma_k^{p + 1}  \\
\\
& \refGE{MkMin} &
\Omega_p(f, x_k; \bar x_{k + 1}) + S_k(\bar x_{k + 1})
- (c + \Delta^{(p)}_{\dom \psi}(f))\gamma_k^{p + 1} \\
\\
& \refGE{def-Delta} &
f(\bar x_{k + 1}) + S_k(\bar x_{k + 1}) - (c + 2\Delta^{(p)}_{\dom \psi}(f))\gamma_k^{p + 1} \\
\\
& = & F_k(\bar x_{k + 1}) - \delta_{k + 1}. \QR
\ea
$$

Thus, we come to the following minimization scheme.
\beq\label{met-CPTM1}
\ba{|c|}
\hline \\
\mbox{\bf Contracting-Point Tensor Method, I}\\
\\
\hline \\
\ba{l}
\mbox{{\bf Initialization.} Choose $x_0 \in \dom \psi$, $c \geq 0$.}\\
\\
\mbox{\bf Iteration $k \geq 0$.}\\
\\
\mbox{1: Choose $\gamma_k \in (0,1]$.}\\
\\
\mbox{2: For some $\xi_{k + 1} \leq c \gamma_k^{p + 1}$, 
	find $\bar{x}_{k+1}$ satisfying
(\ref{MkMin}).
    }\\
\\
\mbox{3: If $F(\bar x_{k + 1}) \leq F(x_k)$, then set
$x_{k + 1} = \bar x_{k + 1}$.
Else choose $x_{k+1} = x_k$.}\\
\ea\\
\\
\hline
\ea
\eeq

For $p = 1$ and $\psi(\cdot)$ being an indicator function
of a compact convex set, this is well-known Frank-Wolfe
algorithm~\cite{frank1956algorithm}. For $p = 2$, this is
Contracting-Domain Newton Method
from~\cite{doikov2020convex}.

Straightforward consequence of our observations is the following
\BT \label{th-Tens1} 
Let $\gamma_k = \frac{p + 1}{k + p + 1}$. Then, for all iterations $\{ x_k \}_{k \geq 1}$
generated by method~\eqref{met-CPTM1}, we have
$$
\ba{rcl}
F(x_k) - F^{*} & \leq & (p + 1)^{p + 1} \cdot (c + 2\Delta^{(p)}_{\dom \psi}) \cdot k^{-p}.
\ea
$$
\ET
\proof
Let us choose $ A_k = k \cdot (k + 1) \cdot \ldots \cdot
(k + p). $ Then, $a_{k + 1} = A_{k + 1} - A_k = \frac{(p +
1)A_{k + 1}}{k + p + 1}$, and
$$
\ba{rcl}
\gamma_k & = &
\frac{a_{k + 1}}{A_{k + 1}} \; = \; \frac{p + 1}{k + p + 1}.
\ea
$$
Combining~\eqref{FConverg1} with
Theorem~\ref{prop-TensorStep1}, we have
$$
\ba{rcl}
F(x_k) - F^{*} & \leq &
\frac{(c + 2\Delta^{(p)}_{\dom \psi} (f))}{A_k} \sum\limits_{i = 1}^k \frac{a_i^{p + 1}}{A_i^p},
\qquad k \geq 1.
\ea
$$
Since
$$
\ba{rcl}
\frac{1}{A_k} \sum\limits_{i = 1}^k
\frac{a_{i}^{p + 1}}{A_i^p}
& = &
\frac{1}{A_k}
\sum\limits_{i = 1}^k \frac{(p + 1)^{p + 1} A_i}{(p + i)^{p + 1}}
\;\; \leq \;\;
\frac{(p + 1)^{p + 1} k}{A_k}
\;\; \leq \;\; \frac{(p + 1)^{p + 1}}{k^p},
\ea
$$
we get the required inequality.
\qed

It is important, that the required level of accuracy $\xi_{k + 1}$ 
for solving the subproblem is not static: it is
changing with iterations. Indeed, from the practical perspective,
there is no need to use high accuracy during the first iterations,
but it is natural to improve our precision while approaching to the optimum.
Inexact proximal-type tensor methods with dynamic inner accuracies were
studied in~\cite{doikov2020inexact}.

Let us note that the objective $M_k(y)$ from~\eqref{TensorSub}
is generally nonconvex for $p \geq 3$, and it may be nontrivial to look for its
global minimum. Because of that, we propose an alternative condition
for the next point. It requires just to find 
an inexact \textit{stationary point} of $\Omega_p(f, x_k; y)$. That is
a point $\bar x_{k + 1}$, satisfying
\beq \label{TensorStationary}
\ba{rcl}
\la \nabla \Omega_p(f, x_k; \bar x_{k + 1}), v - \bar v_{k + 1} \ra + \psi(v)
& \geq & \psi(\bar v_{k + 1}) - \frac{1}{\gamma_k} \xi_{k + 1}, \quad v \in \dom \psi, \\
\\
\bar x_{k + 1} & = & (1-\gamma_k) x_k + \gamma_k \bar v_{k+1},
\ea
\eeq
for some tolerance value $\xi_{k + 1} \geq 0$.

\BT \label{TensorStep2} 
Let point $\bar x_{k + 1}$ satisfy condition~\eqref{TensorStationary} 
with 
$$
\ba{rcl}
\xi_{k + 1} & \leq & c \gamma_k^{p + 1},
\ea
$$ 
for some constant $c \geq 0$.
Then it satisfies inexact condition~\eqref{Inex2} of
the Conceptual Contracting-Point Method with 
$$
\ba{rcl}
\delta_{k + 1} & = & (c + \Gamma^{(p)}_{\dom \psi}(f))\gamma_k^{p + 1}.
\ea
$$
\ET
\proof
Indeed, for any $v \in \dom \psi$, we have
$$
\ba{cl}
& \la \nabla f(\bar x_{k + 1}), v - \bar v_{k + 1} \ra + \psi(v) \\
\\
& \,= \quad
\la \nabla \Omega_p (f, x_k; \bar x_{k + 1}), v - \bar v_{k + 1} \ra + \psi(v) 
+ \la \nabla f(\bar x_{k + 1}) 
- \Omega_p(f, x_k; \bar x_{k + 1}), v - \bar v_{k + 1} \ra \\
\\
& \refGE{TensorStationary} \quad
\psi(\bar v_{k + 1}) - c \gamma_k^p
+ \la \nabla f(\bar x_{k + 1}) 
- \Omega_p(f, x_k; \bar x_{k + 1}), v - \bar v_{k + 1} \ra \\
\\
& \refGE{def-GammaV} \quad
\psi(\bar v_{k + 1}) - (c + \Gamma^{(p)}_{\dom \psi}(f))\gamma_k^p
\;\; = \;\; \psi(\bar v_{k + 1}) - \frac{1}{\gamma_k} \delta_{k + 1}. \QR
\ea 
$$

Now, changing inexactness condition~\eqref{MkMin} 
in method~\eqref{met-CPTM1} by condition~\eqref{TensorStationary}, we come to the following algorithm.
\beq\label{met-CPTM2}
\ba{|c|}
\hline \\
\mbox{\bf Contracting-Point Tensor Method, II}\\
\\
\hline \\
\ba{l}
\mbox{{\bf Initialization.} Choose $x_0 \in \dom \psi$, $c \geq 0$.}\\
\\
\mbox{\bf Iteration $k \geq 0$.}\\
\\
\mbox{1: Choose $\gamma_k \in (0,1]$.}\\
\\
\mbox{2: For some $\xi_{k + 1} \leq c \gamma_k^{p + 1}$, 
	find $\bar{x}_{k+1}$ satisfying
	(\ref{TensorStationary}).
}\\
\\
\mbox{3: If $F(\bar x_{k + 1}) \leq F(x_k)$, then set
	$x_{k + 1} = \bar x_{k + 1}$.
	Else choose $x_{k+1} = x_k$.}\\
\ea\\
\\
\hline
\ea
\eeq

Its convergence analysis is straightforward.
\BT \label{th-Tens2}
Let $A_k \Def k \cdot (k + 1) \cdot \ldots \cdot (k + p)$,
and consequently $\gamma_k = \frac{p + 1}{k + p + 1}$.
Then, for all iterations $\{x_k \}_{k \geq 1}$ of method~\eqref{met-CPTM2}, we have
$$
\ba{rcl}
F(x_k) - F^{*} & \leq &
\ell_k \;\; \leq \;\; 
(p + 1)^{p + 1} \cdot (c + \Gamma^{(p)}_{\dom \psi}(f)) \cdot k^{-p}.
\ea
$$
\ET
\proof
Combining inequality \eqref{FConverg2} with the statement 
of Theorem~\ref{TensorStep2}, we have
$$
\ba{rcl}
F(x_k) - F^{*} & \leq & \ell_k \;\; \leq \;\;
\frac{c + \Gamma^{(p)}_{\dom \psi} (f)}{A_k} 
\sum\limits_{i = 1}^k \frac{a_i^{p + 1}}{A_i^p},
\qquad k \geq 1.
\ea
$$
It remains to use the same reasoning, as in the proof of Theorem~\ref{th-Tens1}.
\qed

\newpage

\section{Inexact contracting Newton method}
\label{sc-Newton}

In this section, let us present
an implementation of our method~\eqref{met-CPTM1} 
for $p = 2$, when at each step we solve the subproblem inexactly
by a variant of first-order Conditional Gradient Method.
The entire algorithm looks as follows.

\beq\label{met-CPHF}
\ba{|c|}
\hline \\
\mbox{\bf Inexact Contracting Newton Method}\\
\\
\hline \\
\ba{l}
\mbox{{\bf Initialization.} Choose $x_0 \in \dom \psi$, $c > 0$.}\\
\\
\mbox{\bf Iteration $k \geq 0$.}\\
\\
\mbox{1: Choose $\gamma_k \in (0,1]$.}\\
\\
\mbox{2: Denote 
	$g_k(v) = \la \nabla f(x_k), v - x_k \ra 
	+ \frac{\gamma_k}{2}\la \nabla^2 f(x_k)(v - x_k), v - x_k \ra$.
}\\
\\
\mbox{3: Initialize inner method $t = 0$, $z_0 = x_k$, $\phi_0(w) \equiv 0$.
}\\
\\
\mbox{4-a: Set $\alpha_t = \frac{2}{t + 2}$.
}\\
\\
\mbox{4-b: Set $\phi_{t + 1}(w) =
	\alpha_t \bigl[ g_k(z_t) + \la \nabla g_k(z_t), w - z_t\ra + \psi(w) \bigr] + (1 - \alpha_t) \phi_t(w)$.
}\\
\\
\mbox{4-c: Compute 
	$w_{t + 1} \in \Argmin\limits_{w} 
	\phi_{t + 1}(w)$.
}\\[15pt]
\mbox{4-d: Set
	$z_{t + 1} = \alpha_t w_{t + 1} + (1 - \alpha_t) z_t$.
}\\
\\
\mbox{4-e: 
	If $g_k(z_{t + 1}) + \psi(z_{t + 1}) - \phi_{t + 1}(w_{t + 1}) > c \gamma_k^2$,
	then
}\\
\\
\mbox{
	\qquad \qquad Set $t = t + 1$ and go to 4-a, else go to 5.
}\\
\\
\mbox{5: Set $\bar x_{k + 1} = \gamma_k z_{t + 1} + (1 - \gamma_k) x_k$. 
}\\
\\
\mbox{6: If $F(\bar x_{k + 1}) \leq F(x_k)$, then set
	$x_{k + 1} = \bar x_{k + 1}$.
	Else choose $x_{k+1} = x_k$.}\\
\ea\\
\\
\hline
\ea
\eeq

We provide an analysis of the total number of \textit{oracle calls} for $f$ (step 2)
and the total number of \textit{linear minimization oracle calls}
for the composite component $\psi$ (step 4-c), 
required to solve problem~\eqref{MainProblem} 
up to the given accuracy level.

\BT \label{th-CPHF} 
Let $\gamma_k = \frac{3}{k + 3}$. 
Then, for iterations $\{ x_k \}_{k \geq 1}$ generated by method~\eqref{met-CPHF},
we have
\beq \label{eq-CPHF-conv1}
\ba{rcl}
F(x_k) - F^{*} & \leq & 
27 \cdot (c + 2 \Delta^{(2)}_{\dom \psi}) \cdot k^{-2}.
\ea
\eeq
Therefore, for any $\varepsilon > 0$, it is enough to perform
\beq \label{eq-CPHF-conv2}
\ba{rcl}
K & = & 
\biggl\lceil   
\sqrt{\frac{27(c + 2\Delta^{(2)}_{\dom \psi}(f))}{\varepsilon}}
\;
\biggr\rceil
\ea
\eeq
iteration of the method, in order to get $F(x_K) - F^{*} \leq \varepsilon$.
And the total number $N_K$ of linear minimization oracle calls  
during these iterations is bounded as
\beq \label{eq-CPHF-conv3}
\ba{rcl}
N_K & \leq & 
2 \cdot \Bigl( 1 + \frac{2 {\cal V}^{(2)}_{\dom \psi}(f)}{c} \Bigr)
\cdot \Bigl(  1 + \frac{  27(c + 2 \Delta^{(2)}_{\dom\psi}(f))   }{\varepsilon} \Bigr).
\ea
\eeq
\ET
\proof
Let us fix arbitrary iteration $k \geq 0$ of our method and consider the following objective:
$$
\ba{rcl}
m_k(v) & = & g_k(v) + \psi(v) \\
\\
& = & \la \nabla f(x_k), v - x_k \ra 
+ \frac{\gamma_k}{2} \la \nabla^2 f(x_k)(v - x_k), v - x_k \ra
+ \psi(v).
\ea
$$
We need to find the point $\bar v_{k + 1}$ such that
\beq \label{mkBound}
\ba{rcl}
m_k(\bar v_{k + 1}) - m_k^{*} & \leq & c \gamma_k^{2}.
\ea
\eeq
Note that if we set 
$\bar x_{k + 1} := \gamma_k \bar v_{k + 1} + (1 - \gamma_k) x_k$,
then from~\eqref{mkBound} we obtain bound~\eqref{MkMin} satisfied
with $\xi_{k + 1} = c \gamma_k^{3}$.
Thus we would obtain iteration 
of Algorithm~\eqref{met-CPTM1} for $p = 2$, 
and Theorem~\ref{th-Tens1} gives the required rate of convergence~\eqref{eq-CPHF-conv1}.
We are about to show 
that steps 4-a -- 4-e of our algorithm are aiming to find such point $\bar v_{k + 1}$.

Let us introduce auxiliary sequences $A_t \Def t \cdot (t + 1)$ 
and $a_{t + 1} \Def A_{t + 1} - A_t$ for $t \geq 0$. 
Then, $\alpha_t \equiv \frac{a_{t + 1}}{A_{t + 1}}$,
and we have the following representation of 
the \textit{Estimating Functions}, for every $t \geq 0$
$$
\ba{rcl}
\phi_{t + 1}(w) & = & 
\frac{1}{A_{t + 1}} \sum\limits_{i = 0}^{t} a_{i + 1} 
\Bigl[ g_k(z_i) + \la \nabla g_k(z_i), w - z_i \ra + \psi(w) \Bigr].
\ea
$$
By convexity of $g_k(\cdot)$, we have
$$
\ba{rcl}
m_k (w) & \geq & \phi_{t + 1}(w), \qquad w \in \dom \psi.
\ea
$$
Therefore, we obtain the following upper bound for the residual~\eqref{mkBound},
for any $v \in \dom \psi$
\beq \label{SubprEst1}
\ba{rcl}
m_k(v) - m_k^{*}
& \leq & 
m_k(v) - \phi_{t + 1}^{*},
\ea
\eeq
where $\phi_{t + 1}^{*} = \min_{w} \phi_{t + 1}(w) = \phi_{t + 1}(w_{t + 1})$.

Now, let us show by induction, that
\beq \label{SubprInduct}
\ba{rcl}
A_{t}\phi_{t}^{*} & \geq & A_{t} m_k(z_{t})
\; - \;
B_{t}, \qquad t \geq 0,
\ea
\eeq
for $B_t := \frac{\gamma_k {\cal V}_{\dom \psi}^{(2)} (f)}{2} \sum_{i = 0}^t \frac{a_{i + 1}^2}{A_{i + 1}}$.
It obviously holds for $t = 0$. Assume that it holds for some $t \geq 0$. Then,
$$
\ba{rcl}
A_{t + 1} \phi_{t + 1}^{*} & = & A_{t + 1} \phi_{t + 1}(w_{t + 1}) \\
\\
& = & A_t \phi_t(w_{t + 1}) 
+ a_{t + 1} \bigl[  g_k(z_{t}) 
+ \la \nabla g_k(z_{t}), w_{t + 1} - z_{t} \ra + \psi(w_{t + 1}) \bigr] \\
\\
& \refGE{SubprInduct} &
A_t m_k(z_t) 
+ a_{t + 1} \bigl[  g_k(z_{t}) 
+ \la \nabla g_k(z_{t}), w_{t + 1} - z_{t} \ra + \psi(w_{t + 1}) \bigr] - B_t \\
\\
& = &
A_{t + 1} \Bigl[ g_k(z_{t}) + \alpha_t \la \nabla g_k(z_t), w_{t + 1} - z_t \ra 
+ \alpha_t \psi(w_{t + 1}) + (1 - \alpha_t) \psi(z_t) \Bigr] - B_t \\
\\
& \geq &
A_{t + 1} \Bigl[ g_k(z_{t}) + \alpha_t \la \nabla g_k(z_t), w_{t + 1} - z_t \ra 
+ \psi(z_{t + 1}) \Bigr] - B_t.
\ea
$$
Note, that
$$
\ba{rcl}
g_k(z_{t + 1}) & = & g_k(z_t + \alpha_t(w_{t + 1} - z_t)) \\
\\
& = & 
g_k(z_t) + \alpha_t \la \nabla g_k(z_t), w_{t + 1} - z_t \ra
+ \frac{\alpha_t^2 \gamma_k}{2} \la \nabla^2 f(x_k)(w_{t + 1} - z_t), w_{t + 1} - z_t \ra.
\ea
$$
Therefore, we obtain
$$
\ba{rcl}
A_{t + 1} \phi_{t + 1}^{*}
& \geq & A_{t + 1} m_k(z_{t + 1}) - B_t 
- \frac{a_{t + 1}^2}{A_{t}} \cdot \frac{\gamma_k {\cal V}_{\dom \psi}^{(2)}(f)}{2},
\ea
$$
and this is~\eqref{SubprInduct} for the next step. 
Therefore, we have~\eqref{SubprInduct} established for all $t \geq 0$.

Combining~\eqref{SubprEst1} with~\eqref{SubprInduct}, we get the following
guarantee for the inner steps 4-a -- 4-e:
$$
\ba{rcl}
m_k(z_{t + 1}) - m_k^{*}
& \leq & m_k(z_{t + 1}) - \phi_{t + 1}^{*}
\;\; \leq \;\;
\frac{\gamma_k {\cal V}_{\dom \psi}^{(2)}(f)}{2 A_{t + 1}} 
\sum\limits_{i = 0}^{t} \frac{a_{i + 1}^2}{A_{i + 1}} \\
\\
& \leq & \frac{2 \gamma_k {\cal V}_{\dom \psi}^{(2)}(f)}{{t + 1}}.
\ea
$$
Therefore, all iterations of our method is well-defined.
We exit from the inner loop on step 4-e after
\beq \label{tBound}
\ba{rcl}
t & \geq & \frac{2 {\cal V}_{\dom \psi}^{(2)}(f) }{c \gamma_k} - 1
\;\; = \;\;
\frac{2(k + 3) {\cal V}_{\dom \psi}^{(2)}(f) }{3 c} - 1,
\ea
\eeq
and the point $\bar{v}_{k + 1} \equiv z_{t + 1}$ satisfies~\eqref{mkBound}.

Hence, we obtain~\eqref{eq-CPHF-conv1} and~\eqref{eq-CPHF-conv2}. 
The total number of linear minimization oracle calls
can be estimated as follows
$$
\ba{rcl}
N_K  & \refLE{tBound} &  
\sum\limits_{k = 0}^{K - 1} 
\Bigl( 1 + \frac{2(k + 3) {\cal V}_{\dom \psi}^{(2)}(f)}{3c}  \Bigr)
\;\; = \;\;
K \Bigl( 1 + \frac{{\cal V}^{(2)}_{\dom \psi}(f)
}{3c}\bigl(  K + 5 \bigr)  \Bigr) \\
\\
& \leq &
K^2 \Bigl( 1 + \frac{2 {\cal V}^{(2)}_{\dom \psi}(f)}{c} \Bigr)
\;\; \leq \;\;
2 \cdot \Bigl( 1 + \frac{2 {\cal V}^{(2)}_{\dom \psi}(f)}{c} \Bigr)
\cdot \Bigl( 1 +  \frac{  27(c + 2 \Delta^{(2)}_{\dom\psi}(f))   }{\varepsilon} \Bigr). \QF
\ea 
$$

According to the result of Theorem~\ref{th-CPHF}, 
in order to solve problem~\eqref{MainProblem} up to $\varepsilon > 0$
accuracy,
we need to perform ${\cal O}(\frac{1}{\varepsilon})$ total
computations of step 4-c of the method (estimate~\eqref{eq-CPHF-conv3}).
This is the same amount of linear minimization oracle calls, as
required in the classical Frank-Wolfe algorithm~\cite{nesterov2018complexity}.
However, this estimate can be over-pessimistic for our method.
Indeed, it comes
as the product of the worst-case complexity bounds 
for the outer and the inner optimization schemes.
It seems to be very rare to meet with the worst-case instance at
the both levels simultaneously.
Thus, the practical performance of our method can be much better.

At the same time, the total number of gradient and Hessian
computations is only ${\cal O}(\frac{1}{\varepsilon^{1/2}})$
(estimate~\eqref{eq-CPHF-conv2}).
This can lead to a significant acceleration over first-order
Frank-Wolfe algorithm, 
when the gradient computation is a bottleneck
(see our experimental comparison in the next section).

The only parameter which remains to choose in method~\eqref{met-CPHF}, 
is the tolerance constant $c > 0$.
Note that the right hand side of~\eqref{eq-CPHF-conv3} is convex in $c$.
Hence, its approximate minimization provides us with the following choice
$$
\ba{rcl}
c & = & 2 \sqrt{ {\cal V}_{\dom \psi}^{(2)}(f) \, \Delta^{(2)}_{\dom\psi}(f) }.
\ea
$$
In practical applications, we may not know some of these constants.
However, in many cases they are small. Therefore, an appropriate choice of $c$ is a small constant.

Finally, let us discuss effective implementation of our method, 
when the composite part is $\{0, +\infty\}$-indicator of the standard simplex:
\beq \label{SimplexDef}
\ba{rcl}
\dom \psi \;\; = \;\; \S_n & \Def & \Bigl\{ x \in \R^n_{+} \; : \; \sum\limits_{i = 1}^n x^{(i)} = 1  \Bigr\}.
\ea
\eeq
This is an example of a set with a finite number of \textit{atoms}, which
are the standard coordinate vectors in this case:
$$
\ba{rcl}
\S_n & = & \conv\{e_1, \dots, e_n \}.
\ea
$$
See~\cite{jaggi2013revisiting} for more examples of atomic sets
in the context of Frank-Wolfe algorithm.
The maximization of a convex function over such sets can be implemented
very efficiently, since the maximum is always at the corner (one of the atoms).

At iteration $k \geq 0$ of method~\eqref{met-CPHF}, we need to minimize over $\S_n$ the quadratic function
$$
\ba{rcl}
g_k(v) & = & \la \nabla f(x_k), v - x_k \ra + \frac{\gamma_k}{2} \la \nabla^2 f(x_k)(v - x_k), v - x_k \ra,
\ea
$$
whose gradient is
$$
\ba{rcl}
\nabla g_k(v)
= \nabla f(x_k) + \gamma_k \nabla^2 f(x_k)(v - x_k).
\ea
$$
Assume that we keep the vector $\nabla g_k(z_t) \in \R^n$ for the current point $z_t$, $t \geq 0$
 of the inner process,
 as well as its aggregation 
 $$
 \ba{rcl}
 h_t & \Def & \alpha_t \nabla g_k(z_t) + (1 - \alpha_t) h_{t - 1},
 \qquad 
 h_{-1} \;\; \Def \;\; 0 \in \R^n.
 \ea
 $$
 Then, at step 4-c we need to compute a vector 
 $$
 \ba{rcl}
 w_{t + 1} & \in & \Argmin\limits_{w \in \S_n} \la h_t, w \ra
 \;\; = \;\; 
 \conv\Bigl\{ e_j \; : \;  j \in \Argmin\limits_{1 \leq j \leq n} h_t^{(j)}  \Bigr \}.
 \ea
 $$
It is enough to find an index $j$ of a minimal element of $h_t$ and to set $w_{t + 1} := e_j$.
The new gradient is equal to
$$
\ba{rcl}
\nabla g_k(z_{t + 1})
& \overset{\text{Step 4-d}}{=} & 
\nabla g_k(\alpha_t w_{t + 1} +  (1 - \alpha_t) z_t ) \\
\\
& = &
\alpha_t\Bigl( 
\nabla f(x_k) + \gamma_k \nabla^2 f(x_k) (e_j - x_k)  
\Bigr) + (1 - \alpha_t) \nabla g_k(z_t),
\ea
$$
and the function value can be expressed using the gradient as follows
$$
\ba{rcl}
g_k(z_{t + 1}) & = &  \frac{1}{2}\la \nabla f(x_k) + \nabla g_k(z_{t + 1}), z_{t + 1} - x_k \ra.
\ea
$$
The product $\nabla^2 f(x_k) e_j$ is just $j$-th column of 
the matrix. 
Hence, preparing in advance the following objects: 
$\nabla f(x_k) \in \R^n$, $\nabla^2 f(x_k) \in \R^{n \times n}$ 
and the Hessian-vector product $\nabla^2 f(x_k) x_k \in \R^n$,
we are able to perform iteration of the inner loop (steps 4-a -- 4-e)
very efficiently in ${ \cal O}(n)$ arithmetical operations.

\section{Numerical experiments}
\label{sc-Numerical}

\begin{figure}[h!]
	\centering
	\includegraphics[width=0.21\textwidth ]{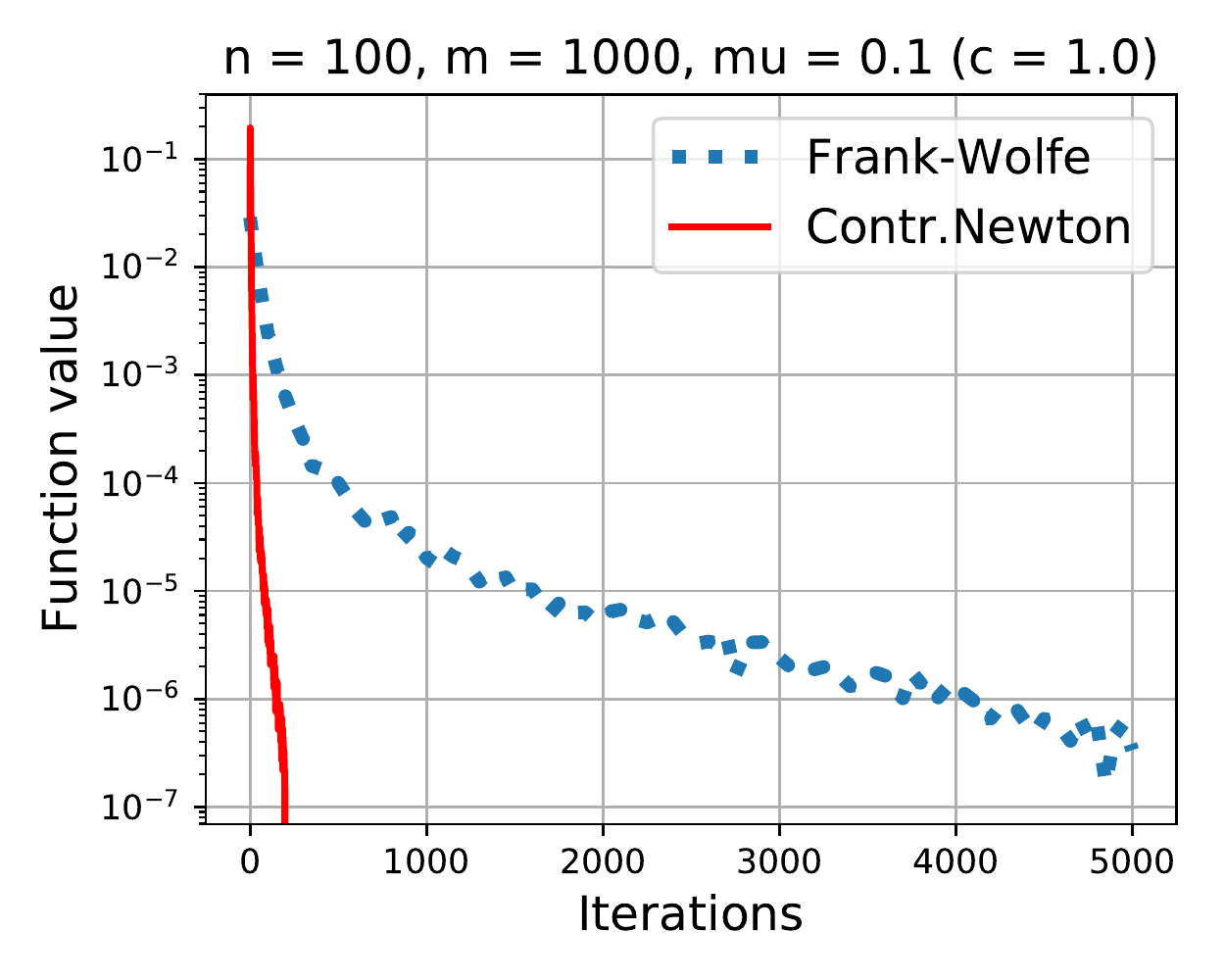}
	\includegraphics[width=0.21\textwidth ]{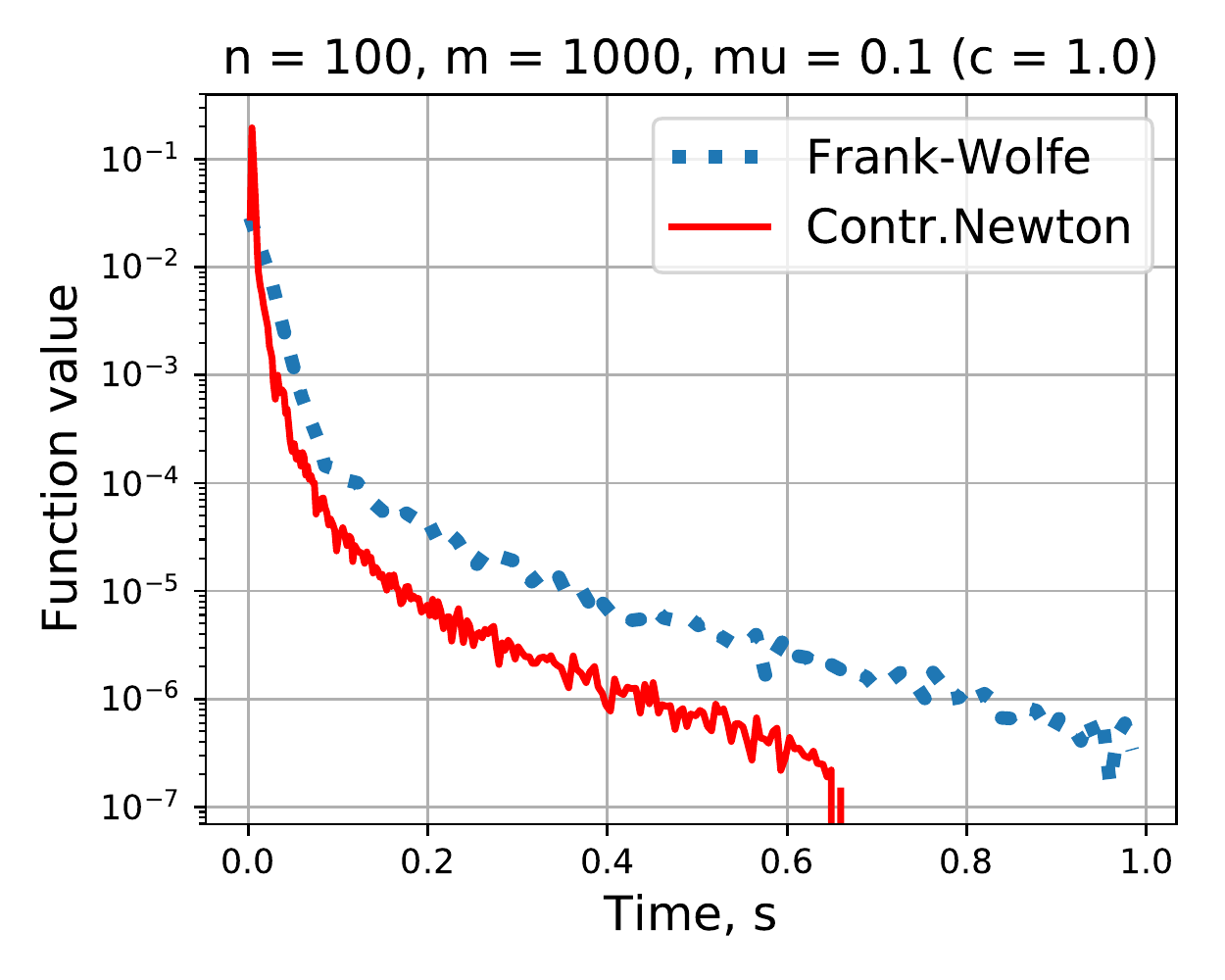}
	\includegraphics[width=0.21\textwidth ]{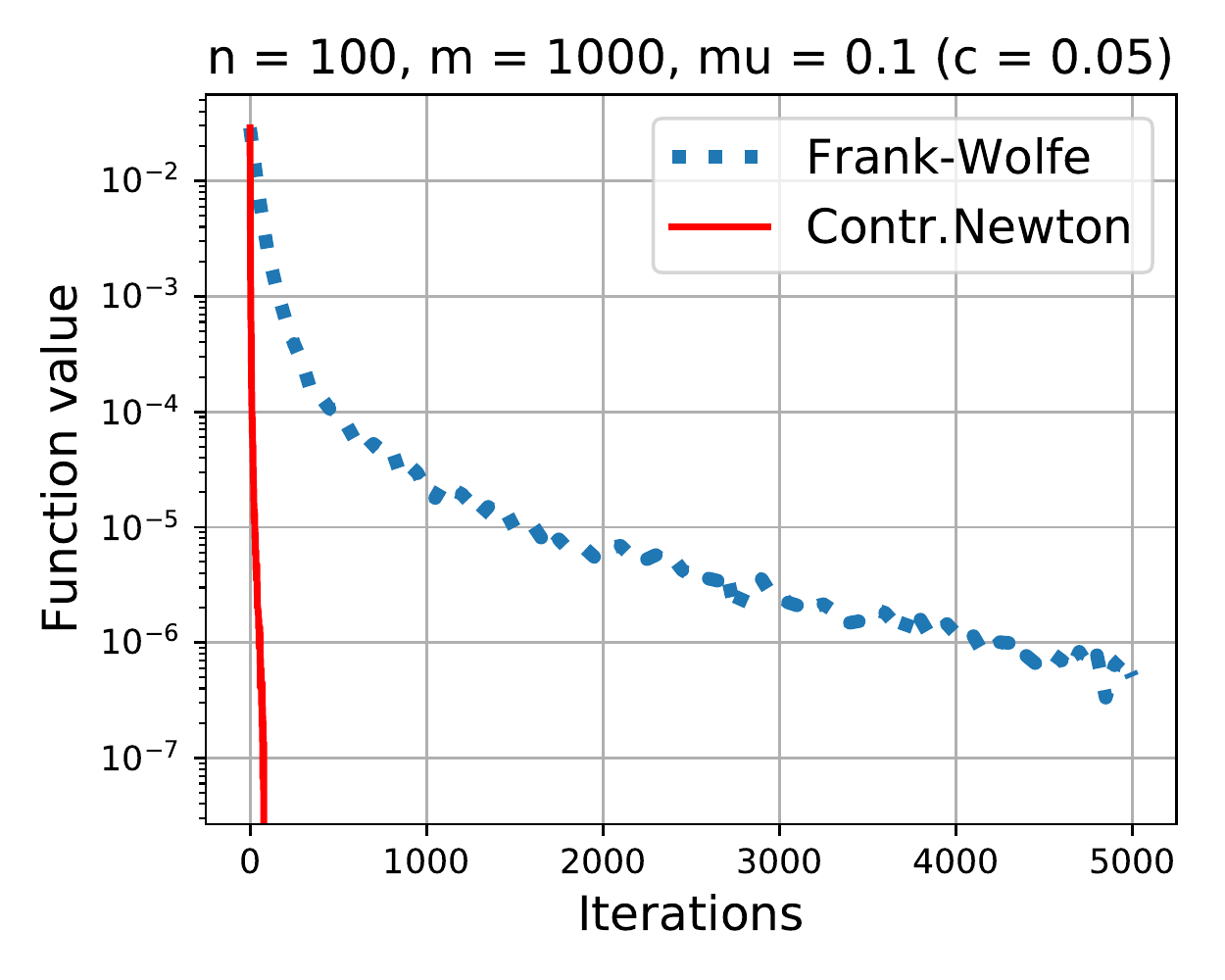}
	\includegraphics[width=0.21\textwidth ]{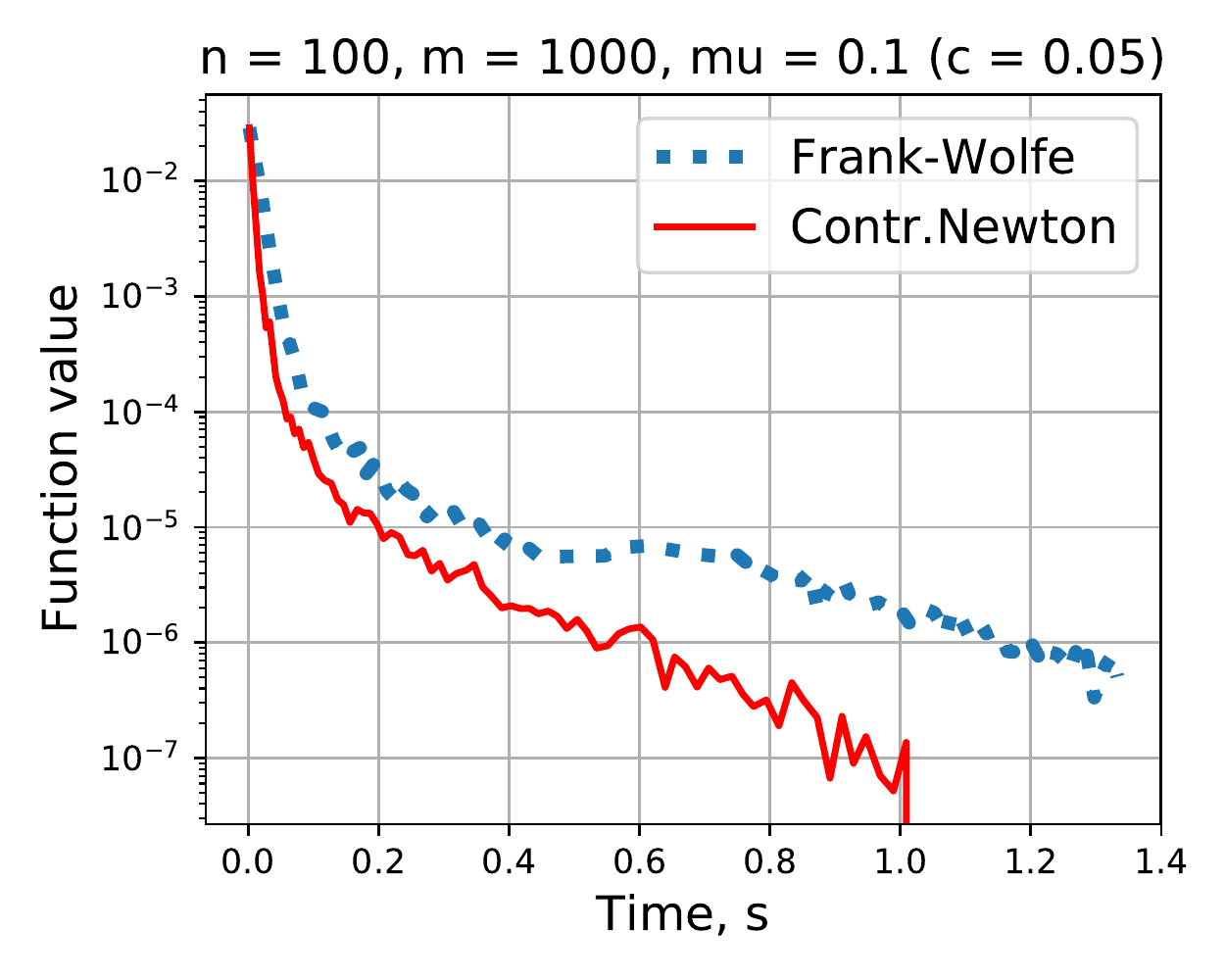}
	
	\includegraphics[width=0.21\textwidth ]{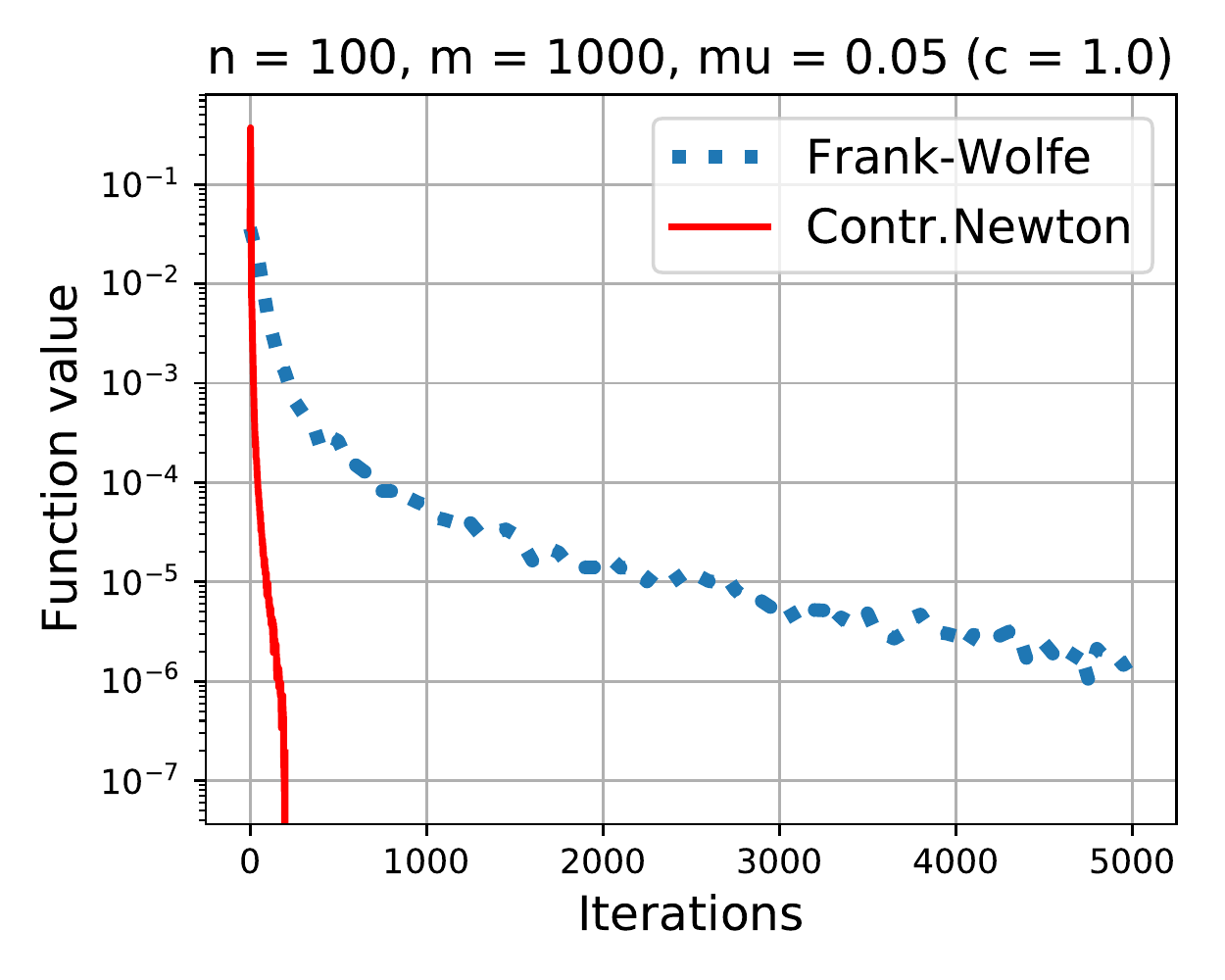}
	\includegraphics[width=0.21\textwidth ]{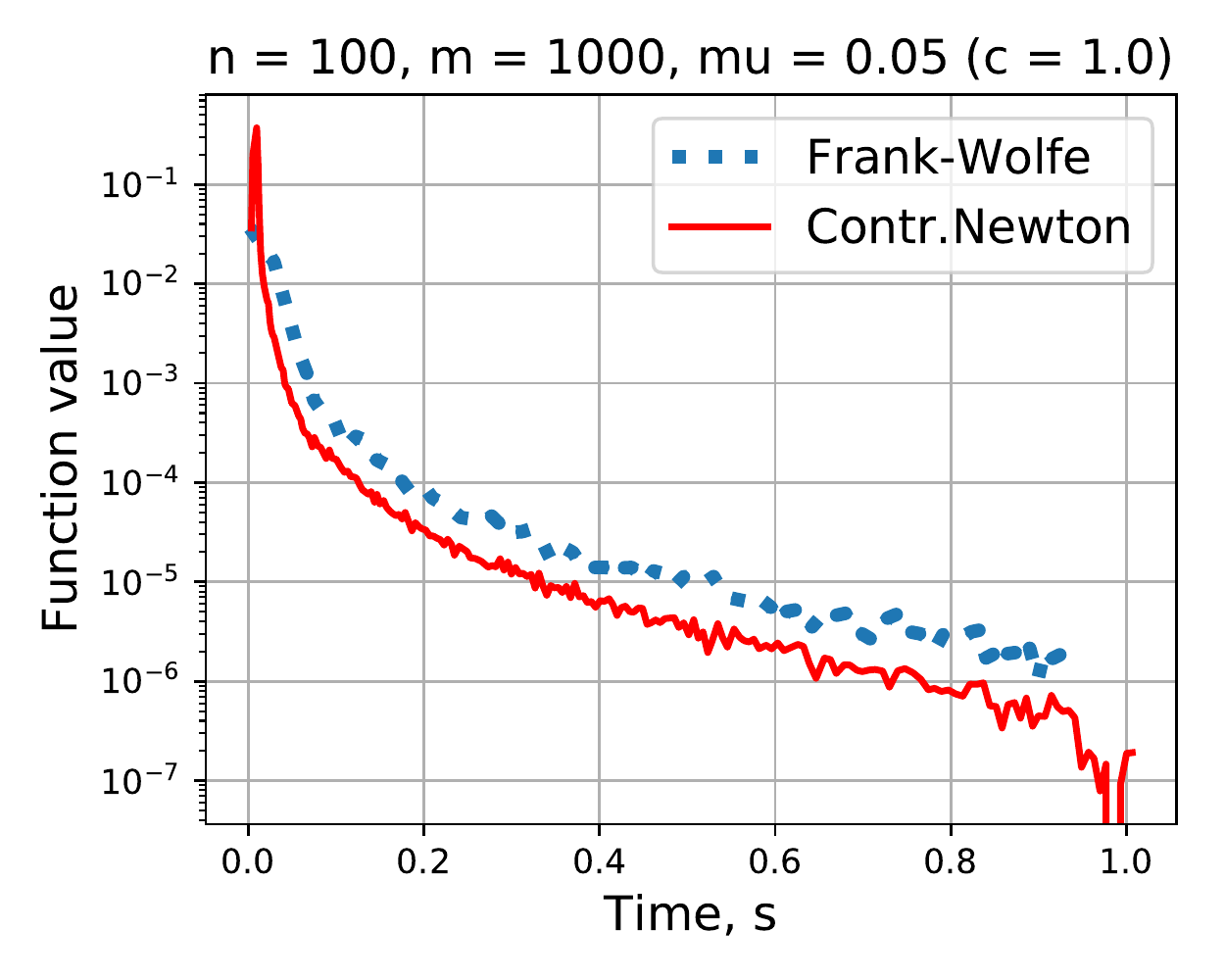}
	\includegraphics[width=0.21\textwidth ]{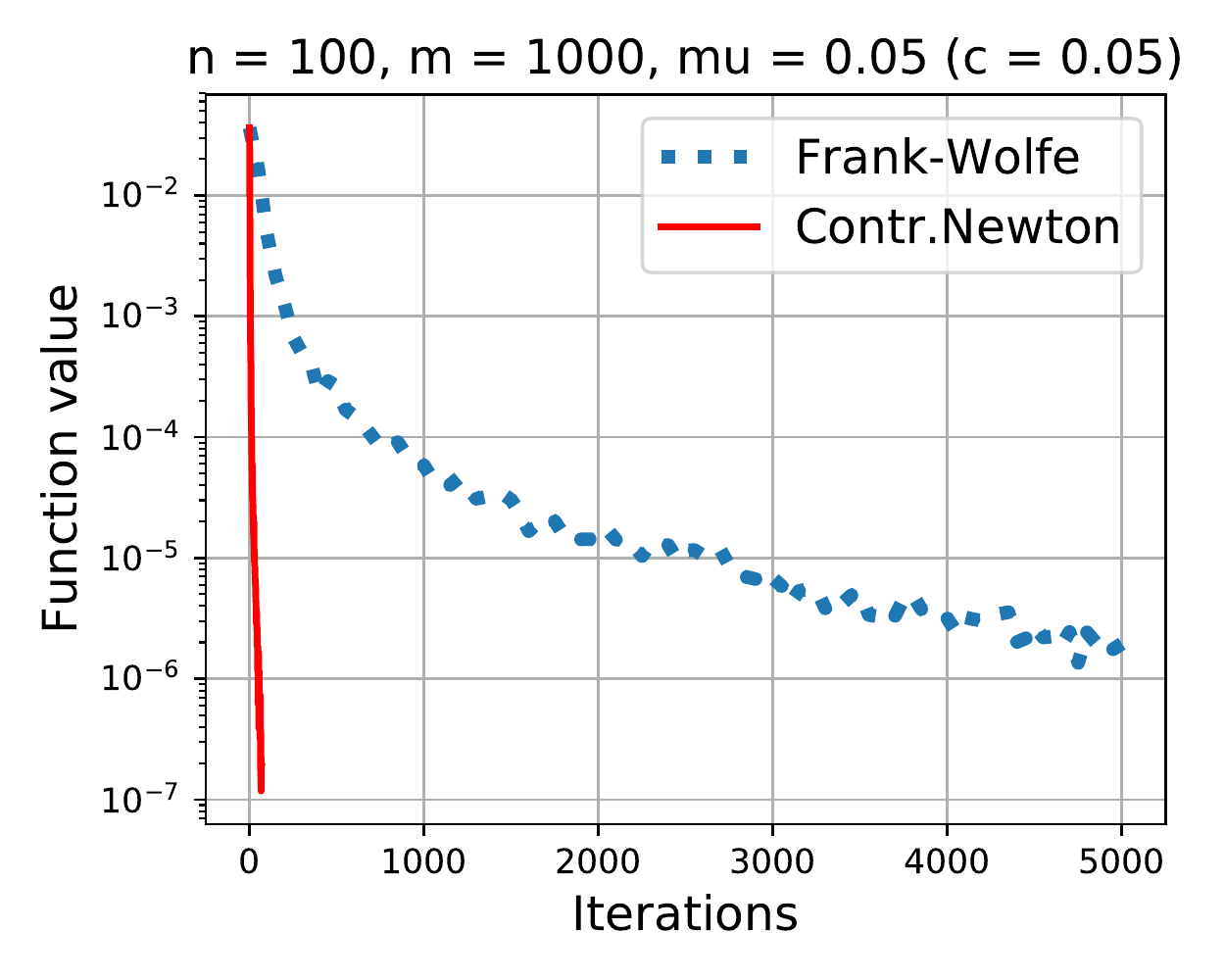}
	\includegraphics[width=0.21\textwidth ]{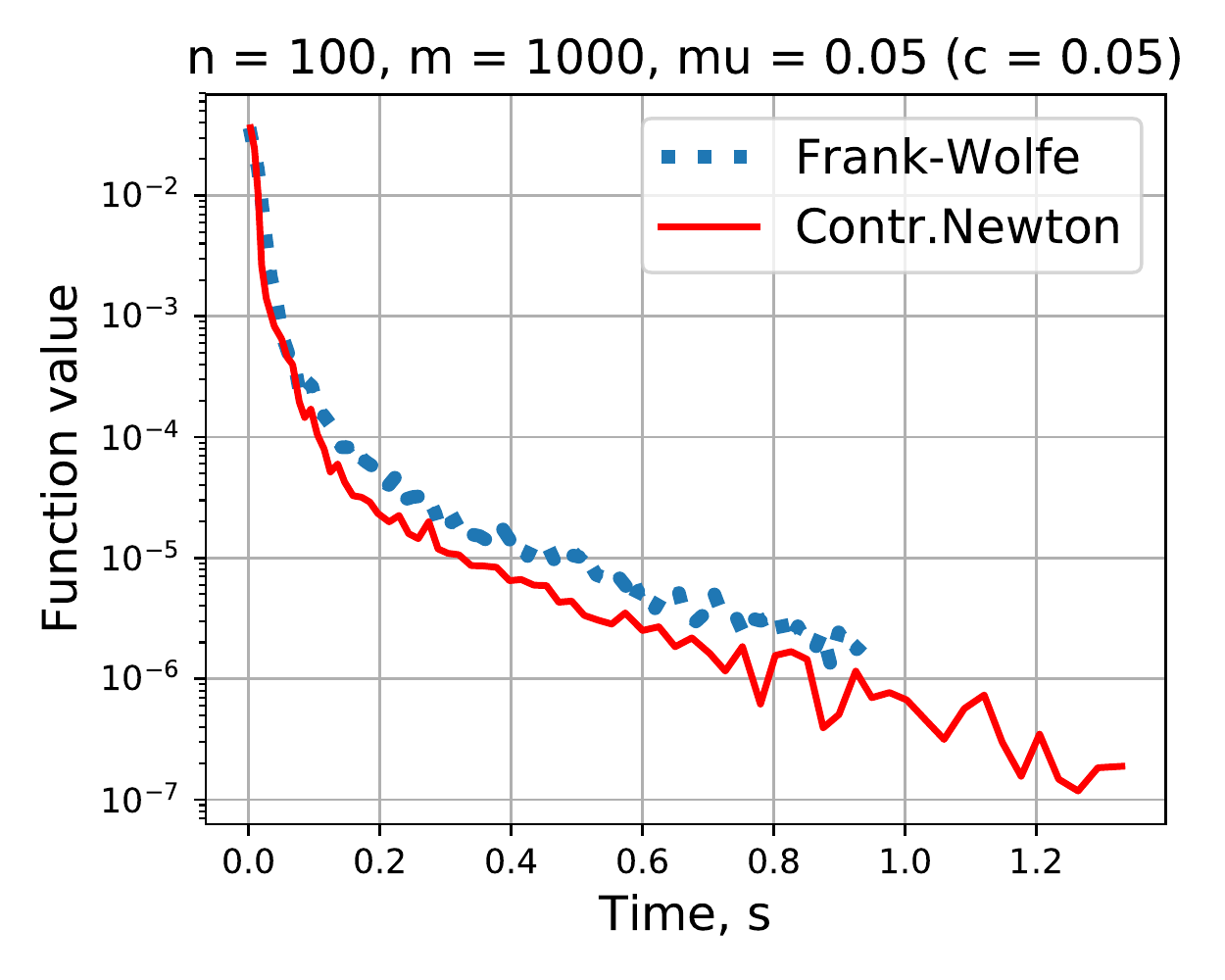}
	
	\caption{$n = 100, \; m = 1000$.}
	\label{Fig1}
\end{figure}

\begin{figure}[h!]
	\centering
	\includegraphics[width=0.21\textwidth ]{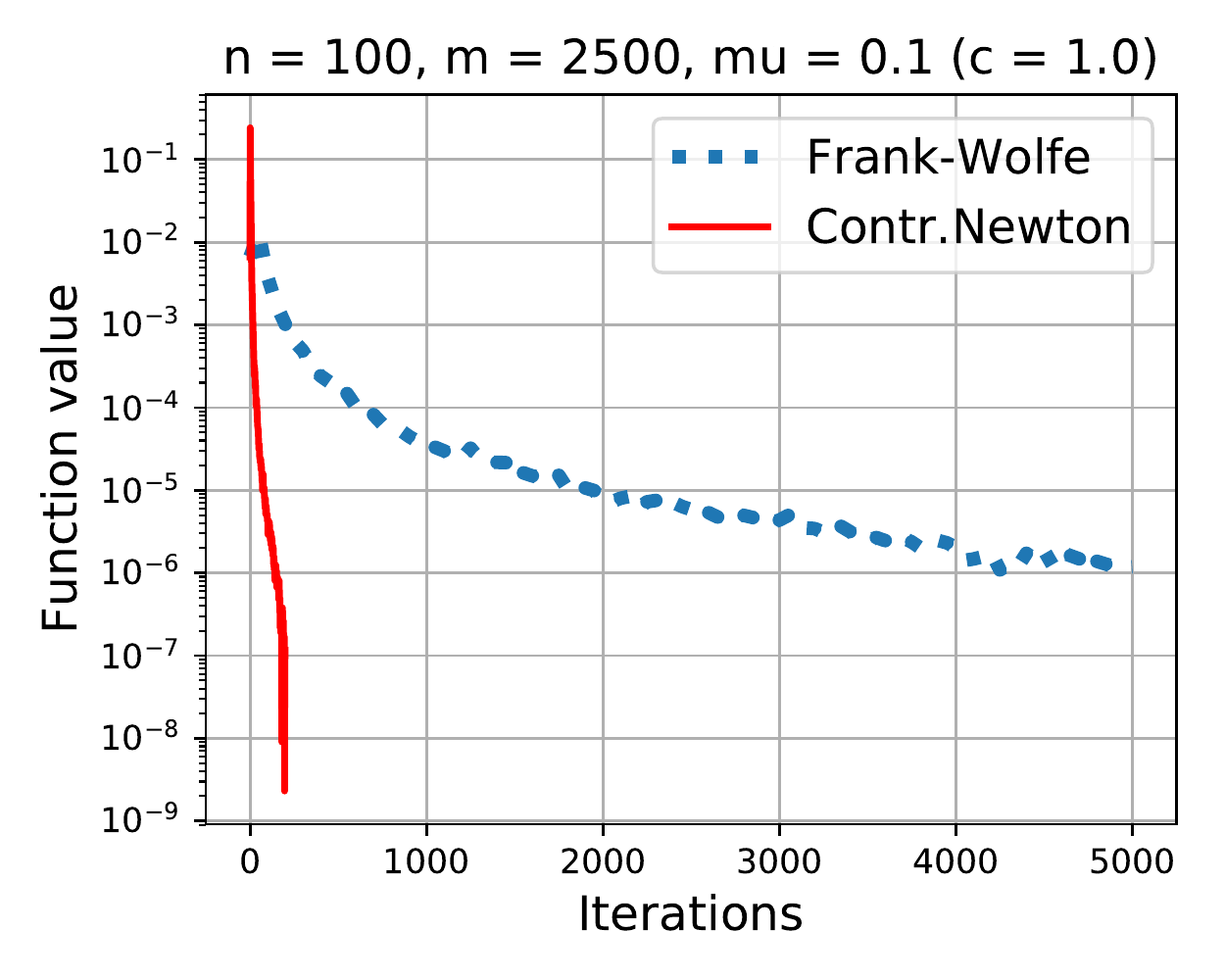}
	\includegraphics[width=0.21\textwidth ]{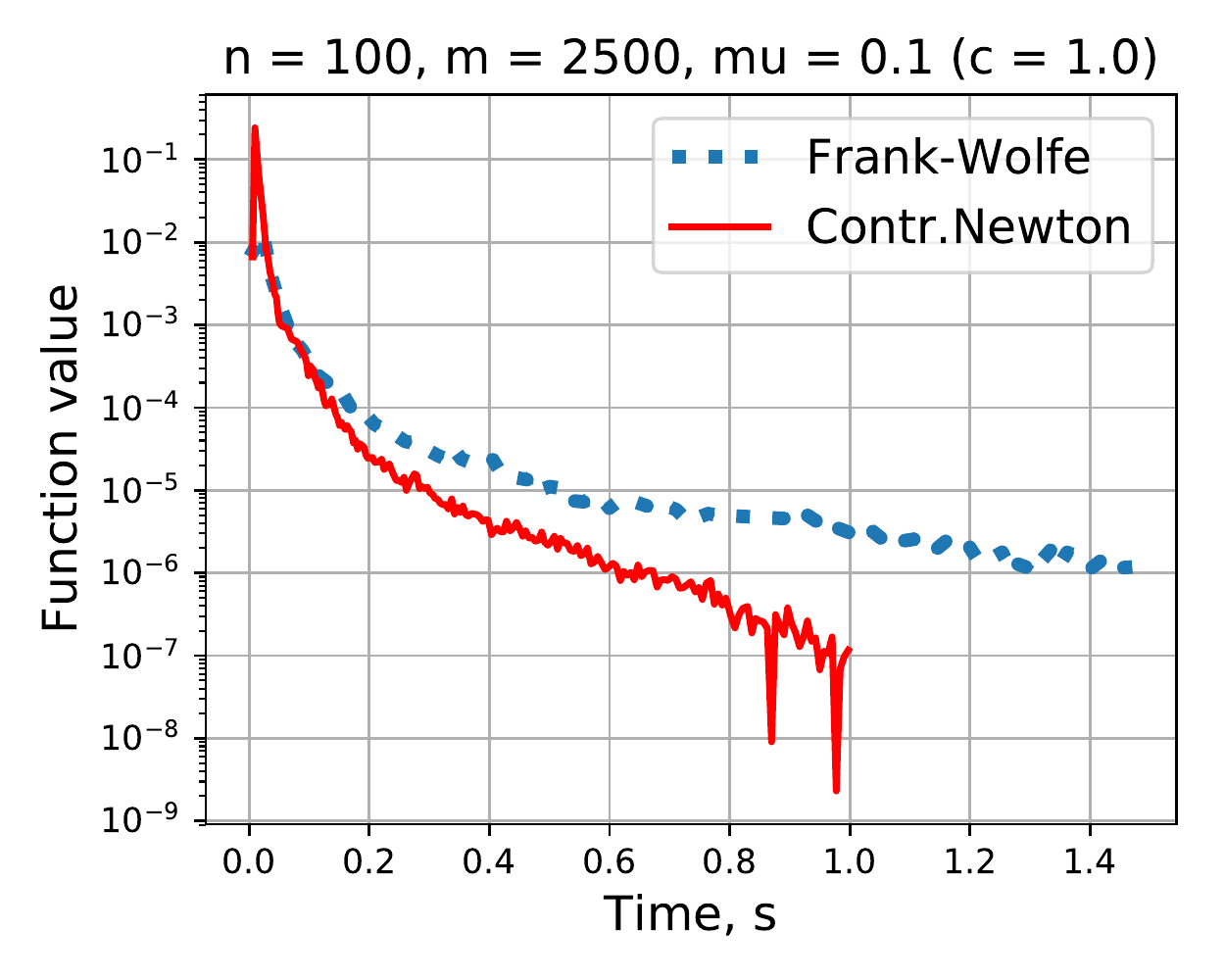}
	\includegraphics[width=0.21\textwidth ]{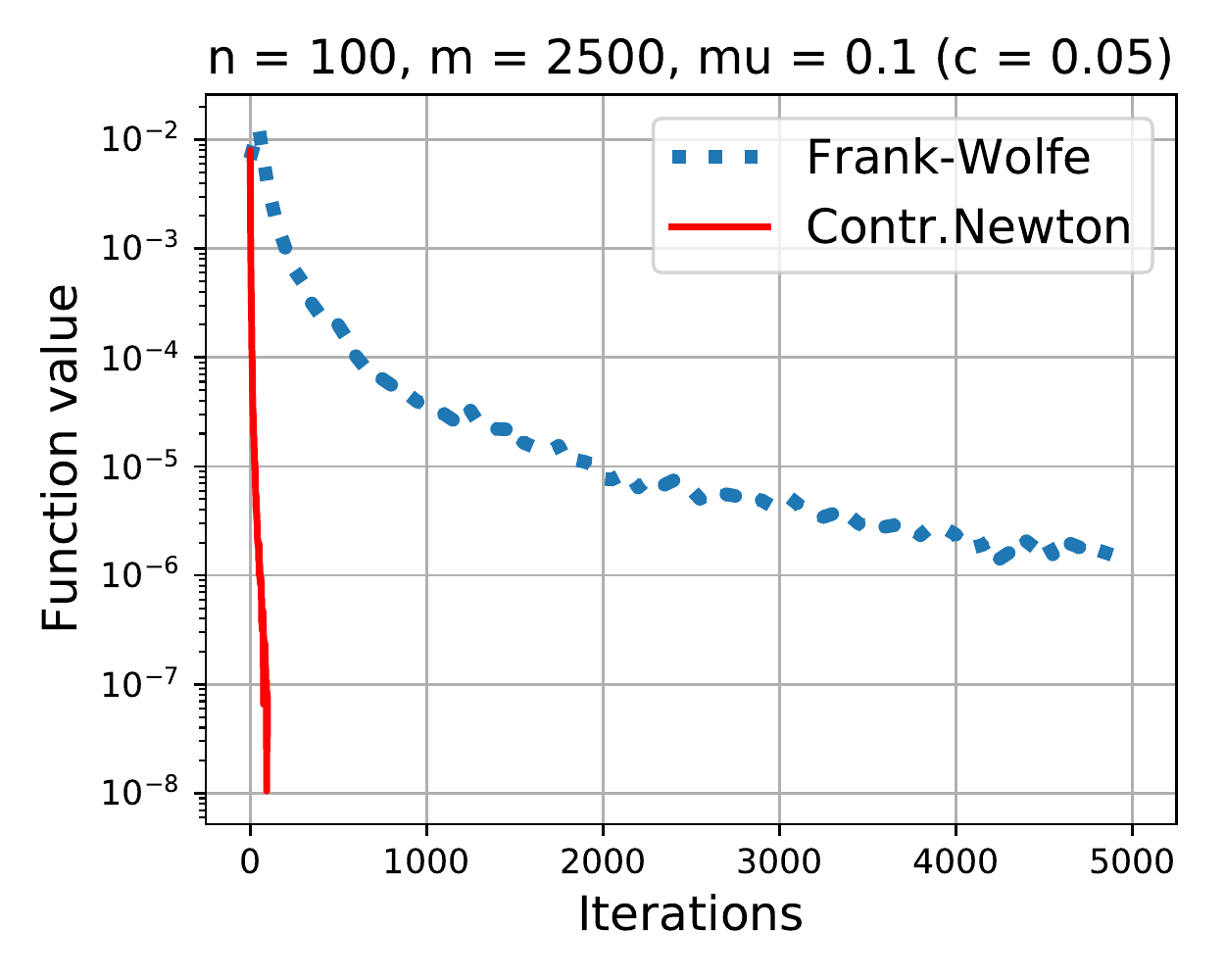}
	\includegraphics[width=0.21\textwidth ]{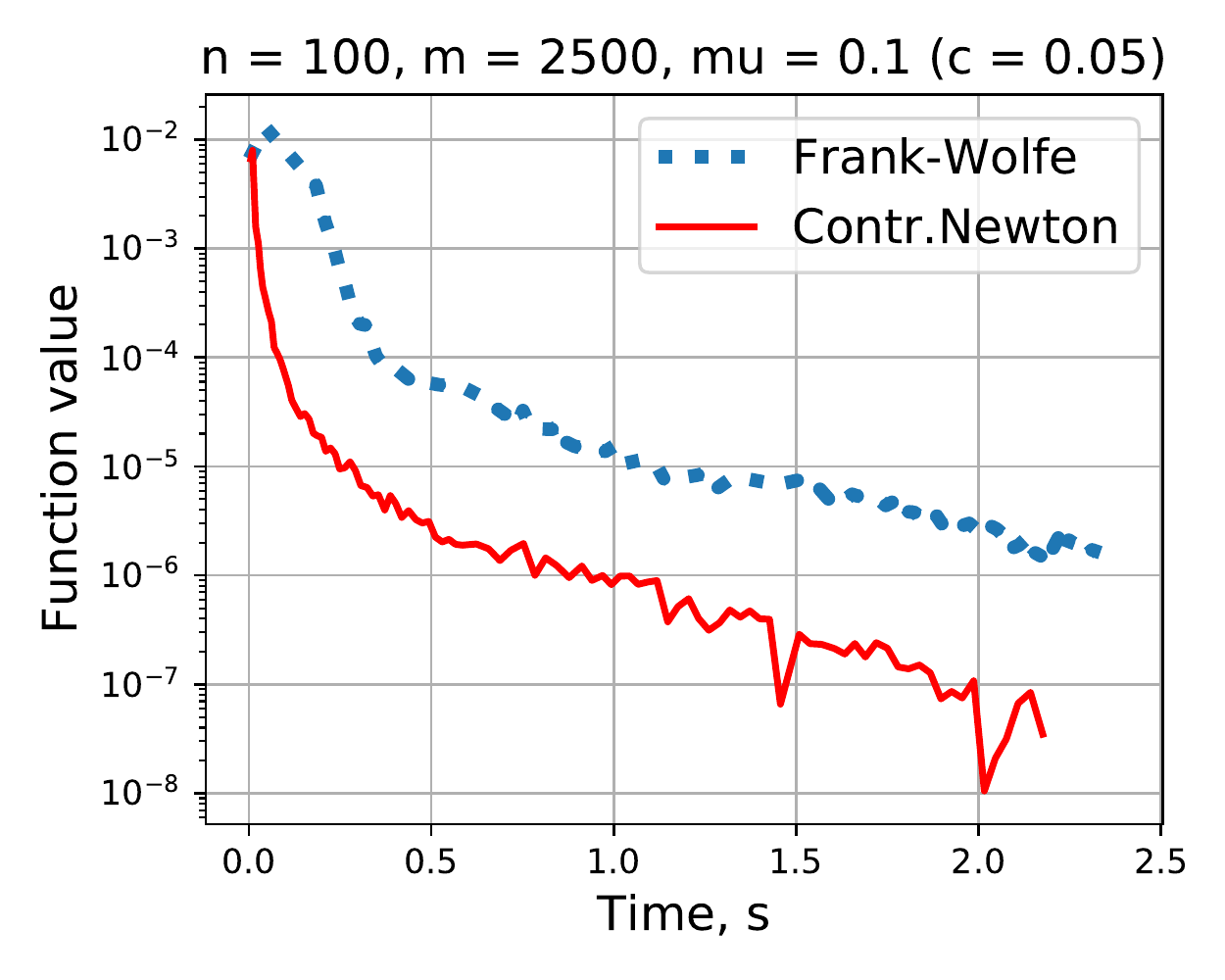}
	
	\includegraphics[width=0.21\textwidth ]{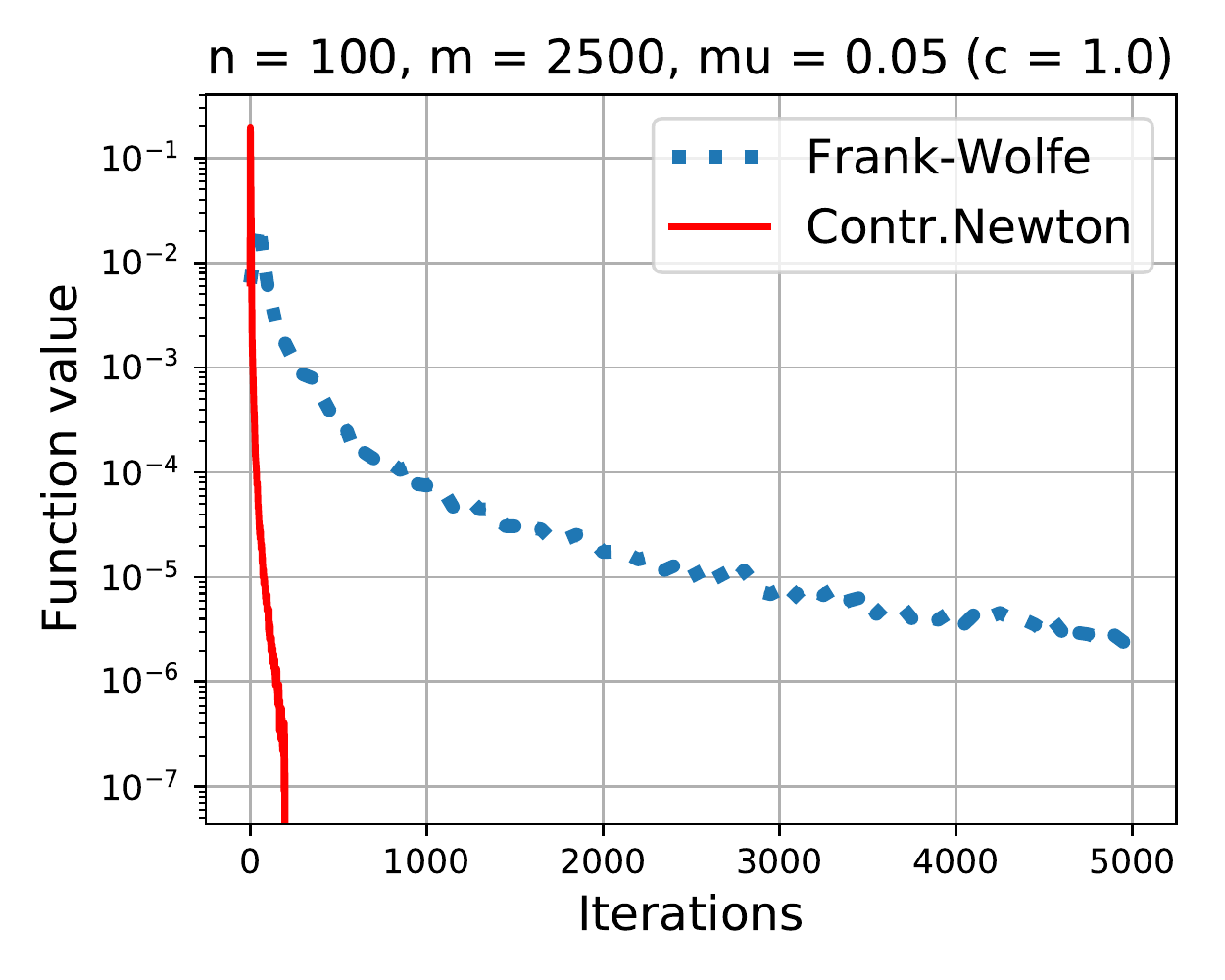}
	\includegraphics[width=0.21\textwidth ]{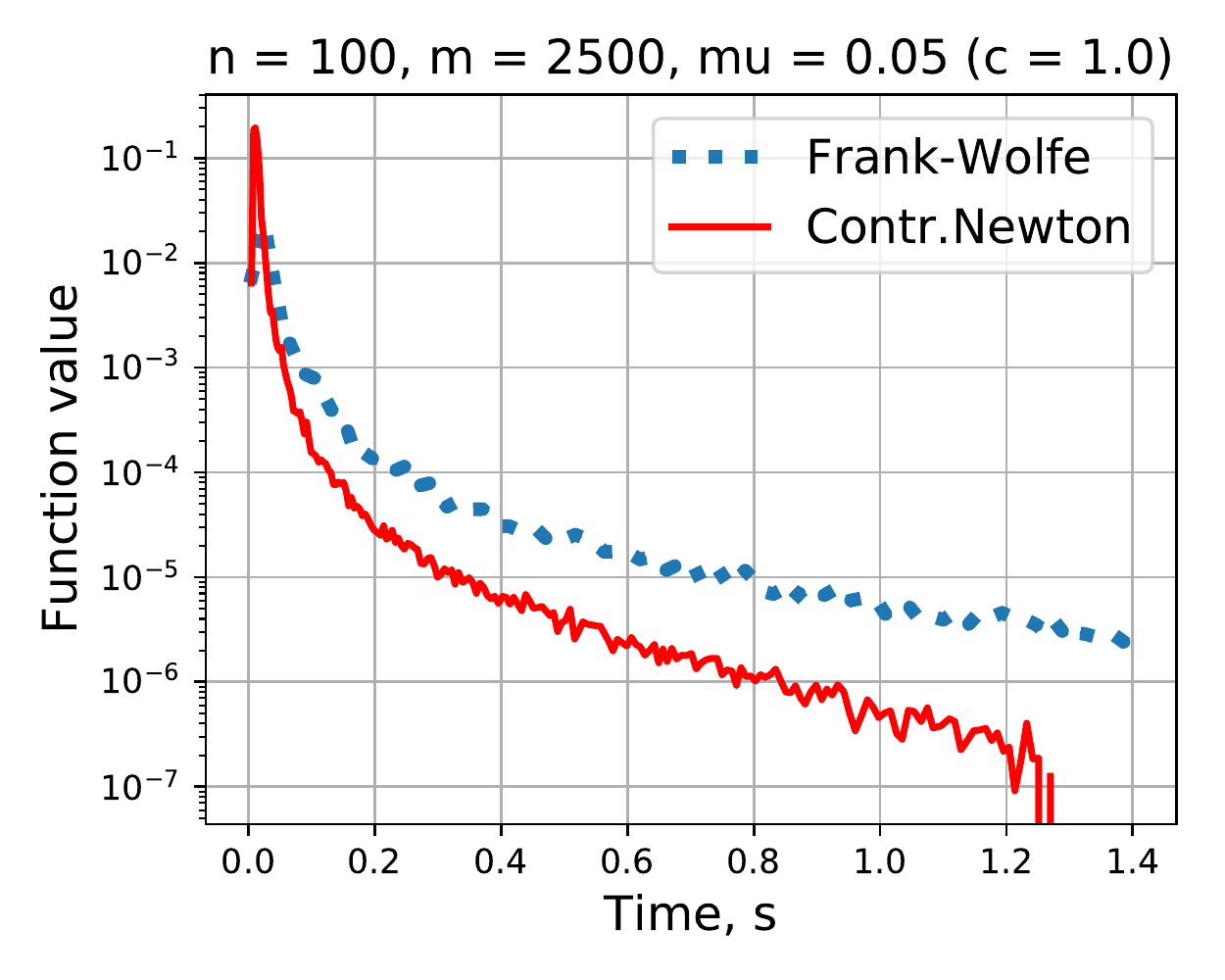}
	\includegraphics[width=0.21\textwidth ]{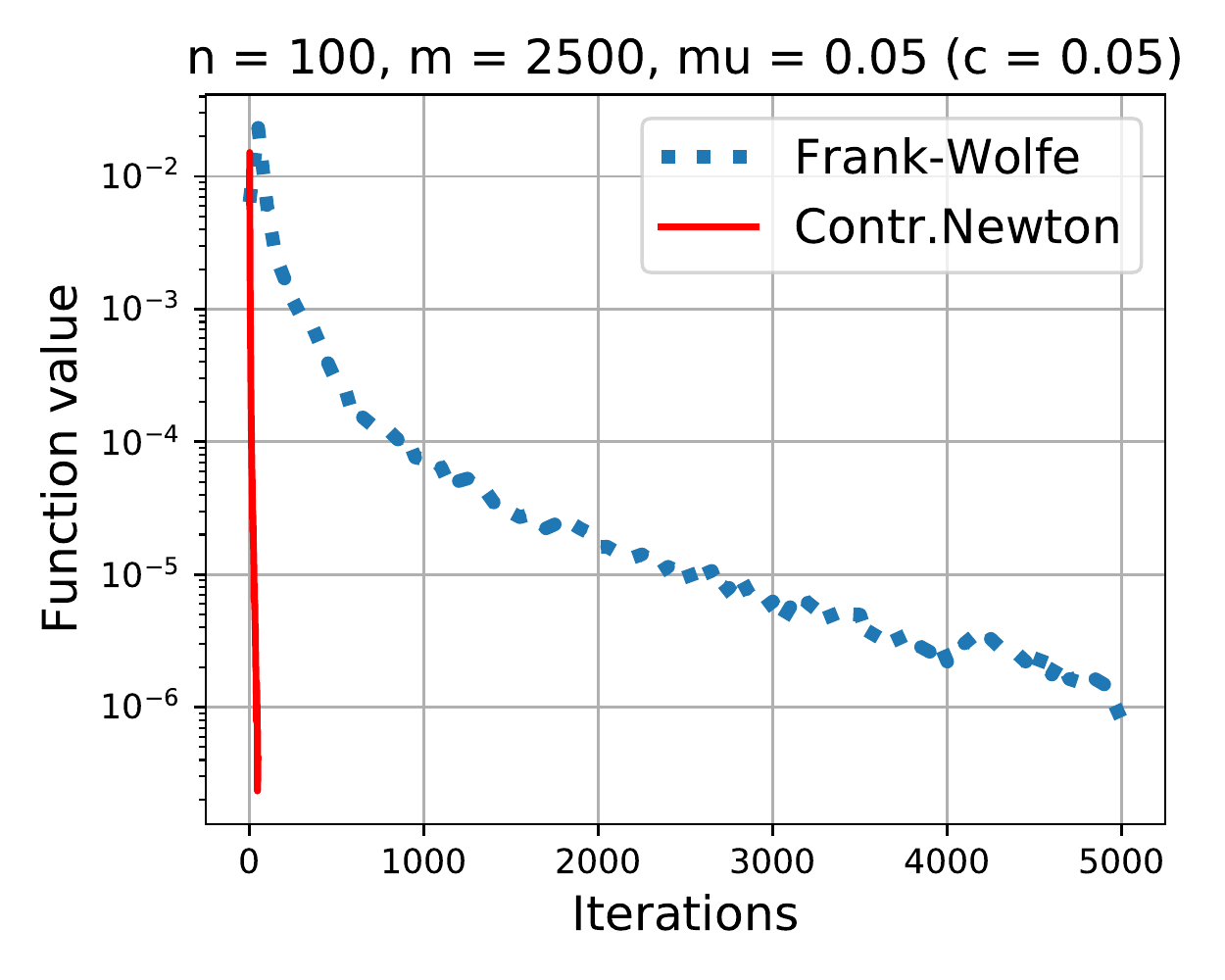}
	\includegraphics[width=0.21\textwidth ]{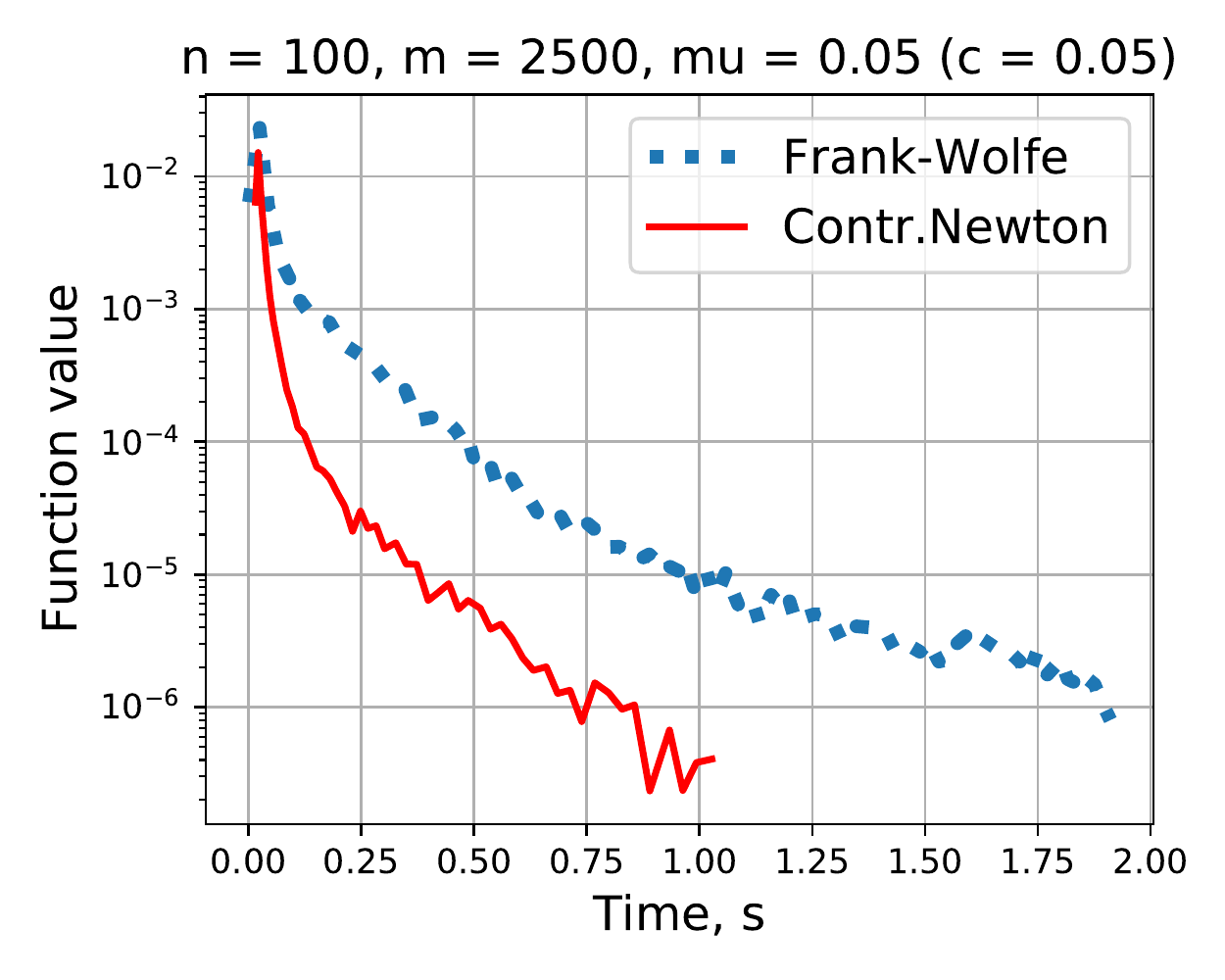}
	
	\caption{$n = 100, \; m = 2500$.}
	\label{Fig2}
\end{figure}

\begin{figure}[h!]
	\centering
	\includegraphics[width=0.21\textwidth ]{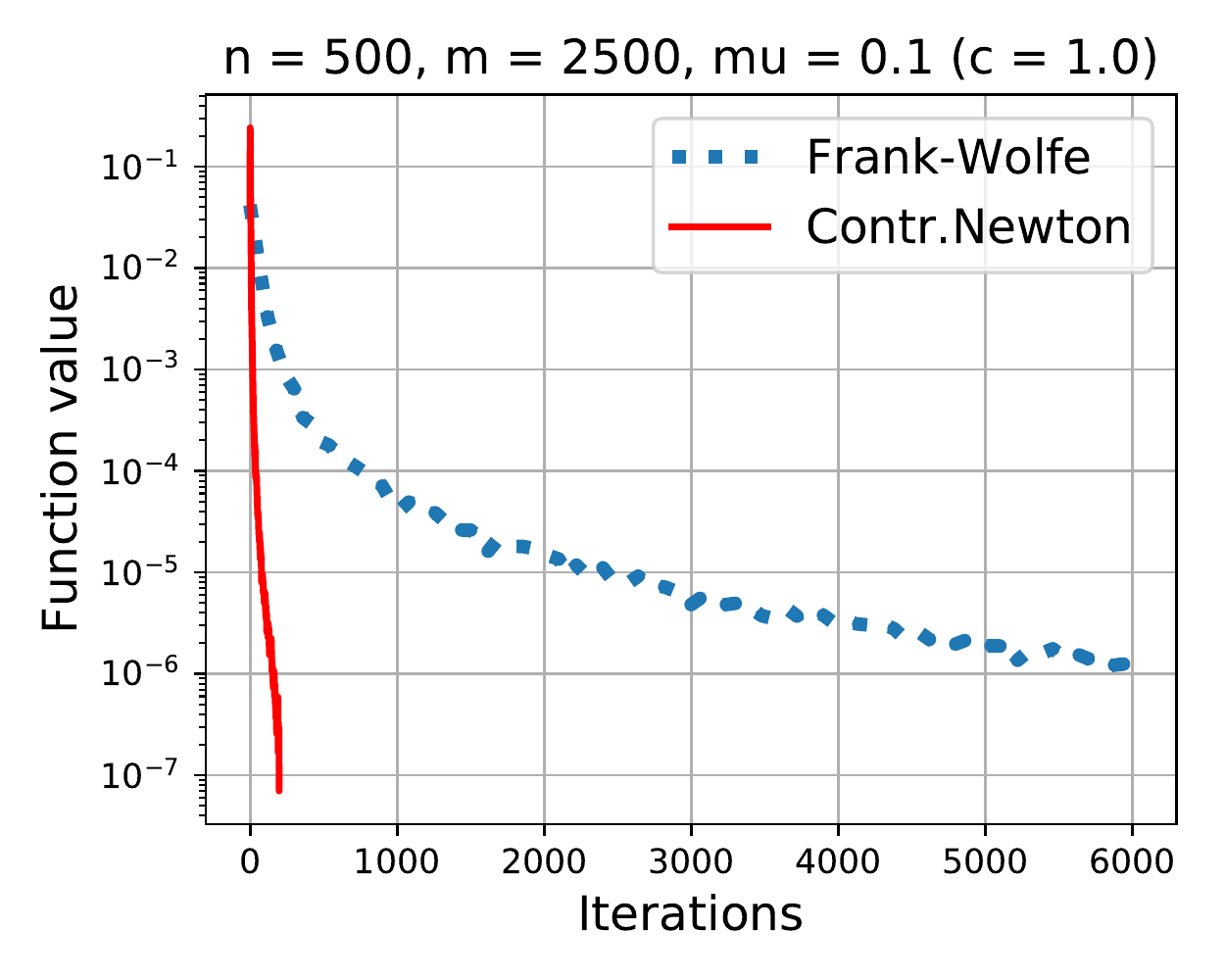}
	\includegraphics[width=0.21\textwidth ]{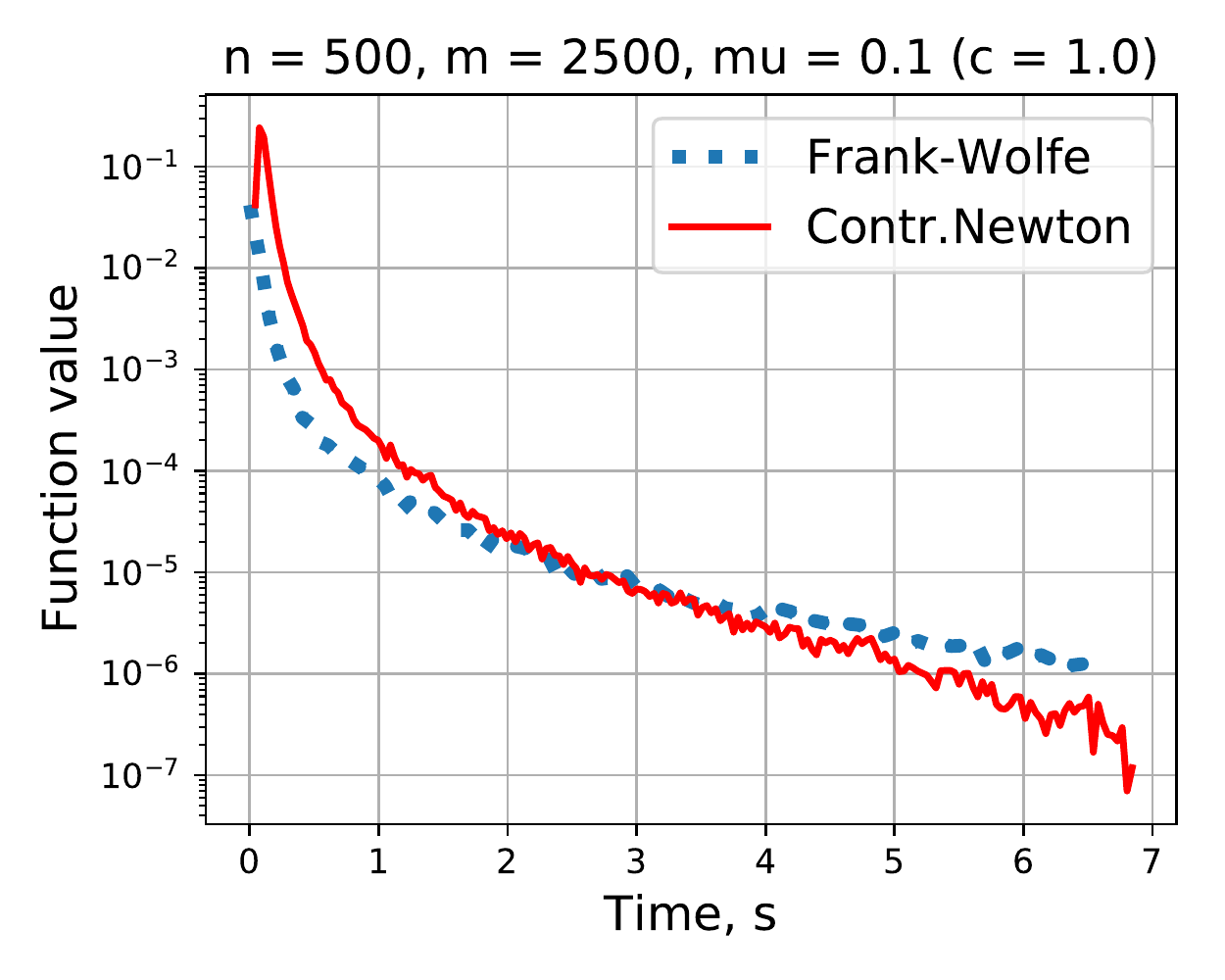}
	\includegraphics[width=0.21\textwidth ]{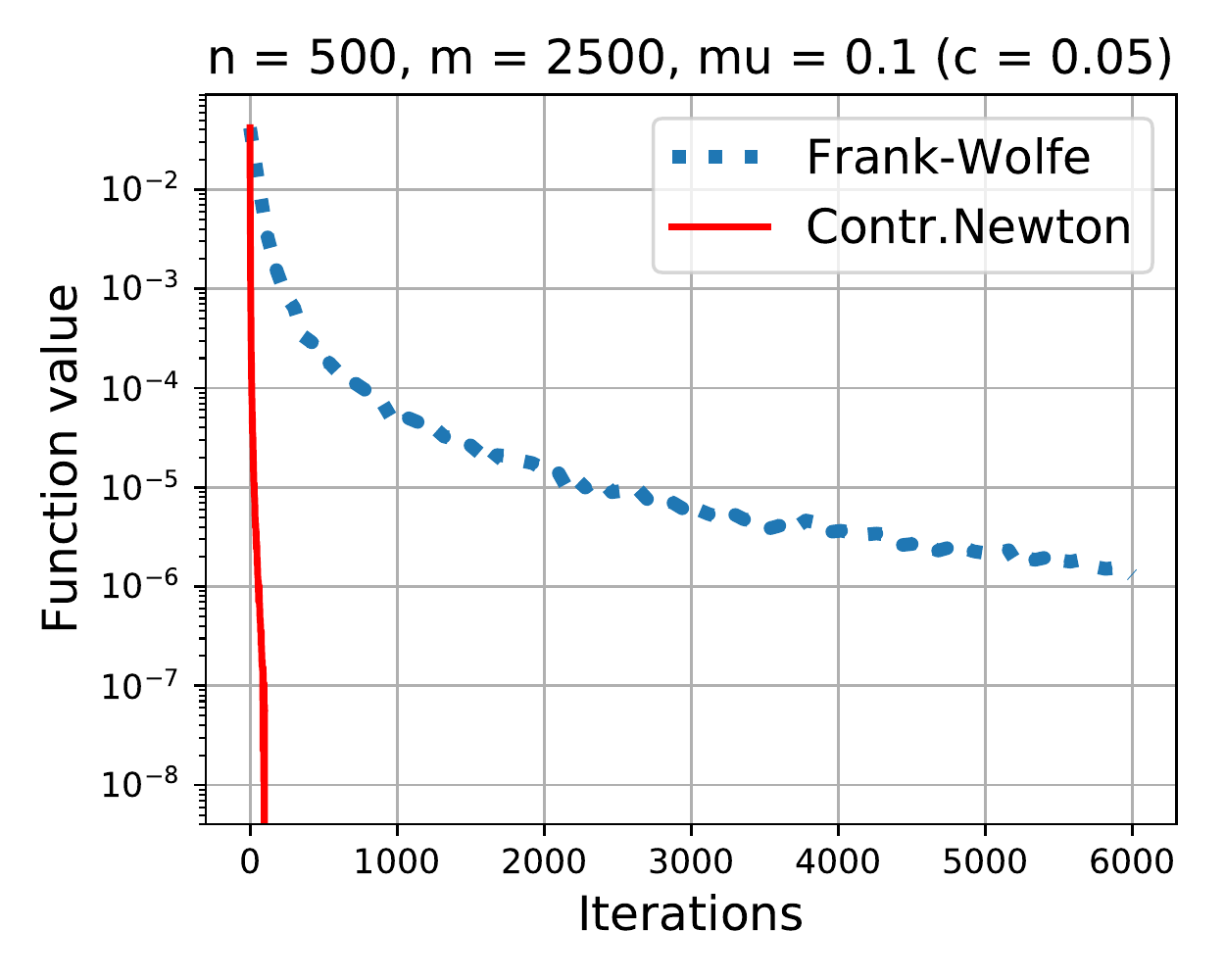}
	\includegraphics[width=0.21\textwidth ]{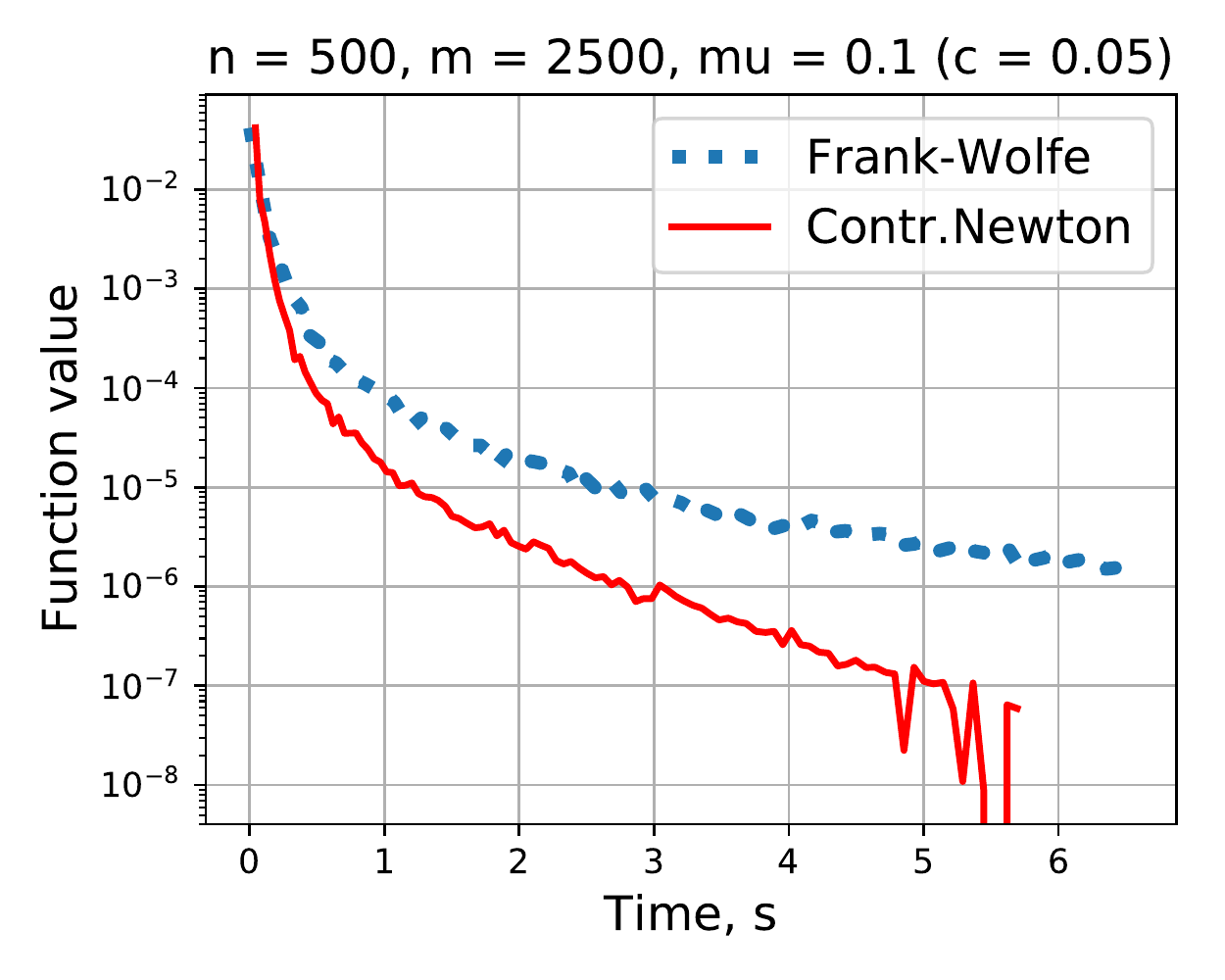}
	
	\includegraphics[width=0.21\textwidth ]{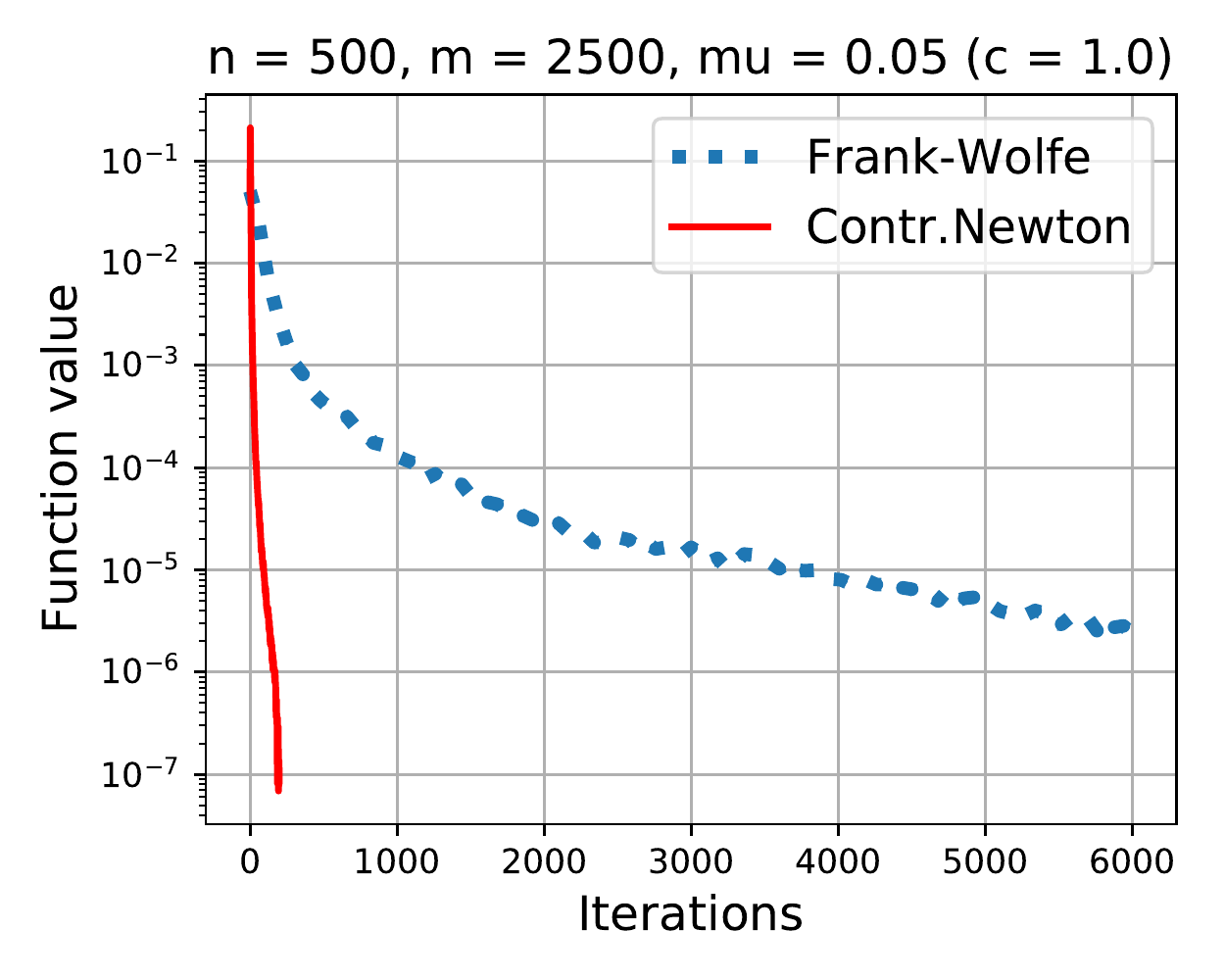}
	\includegraphics[width=0.21\textwidth ]{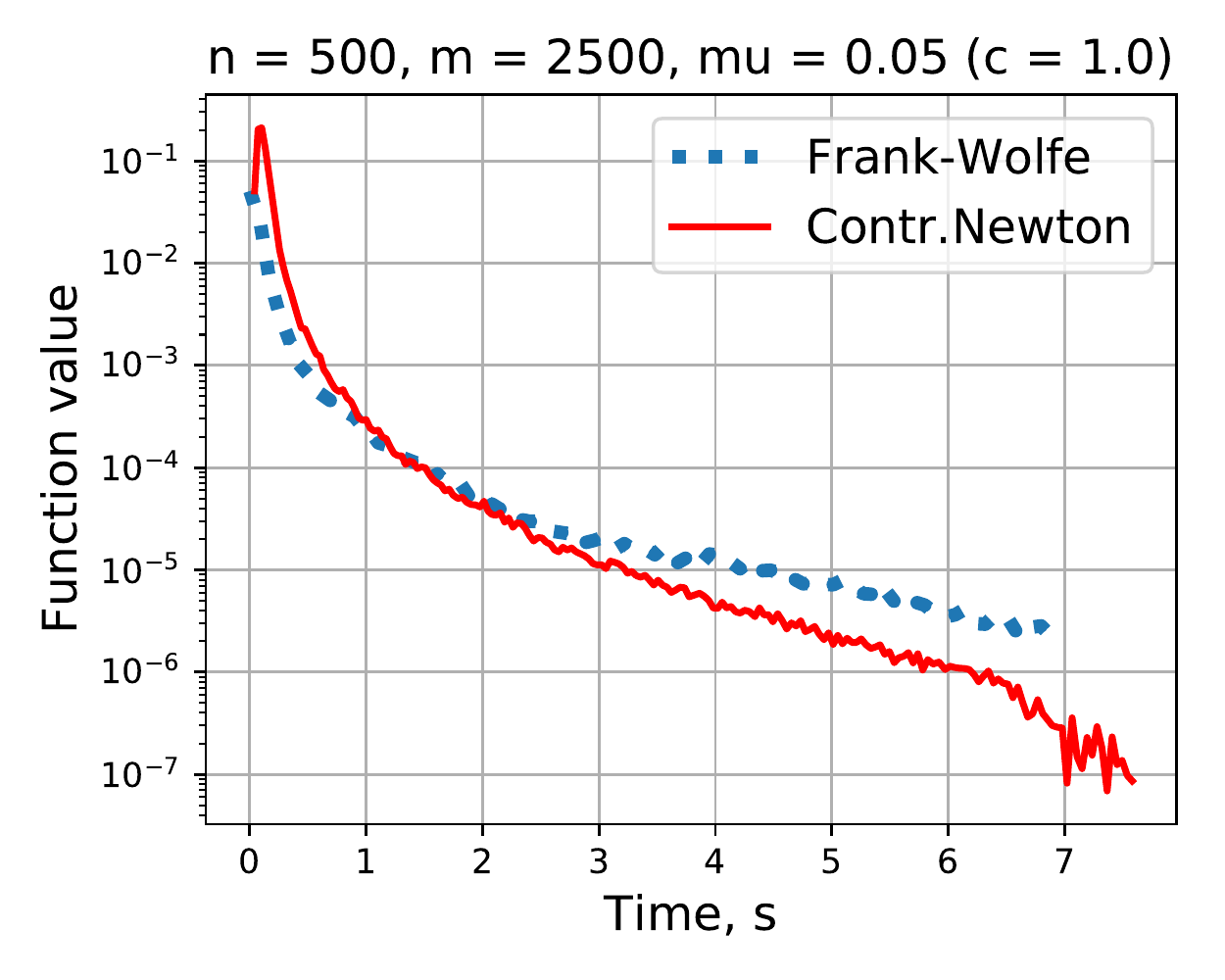}
	\includegraphics[width=0.21\textwidth ]{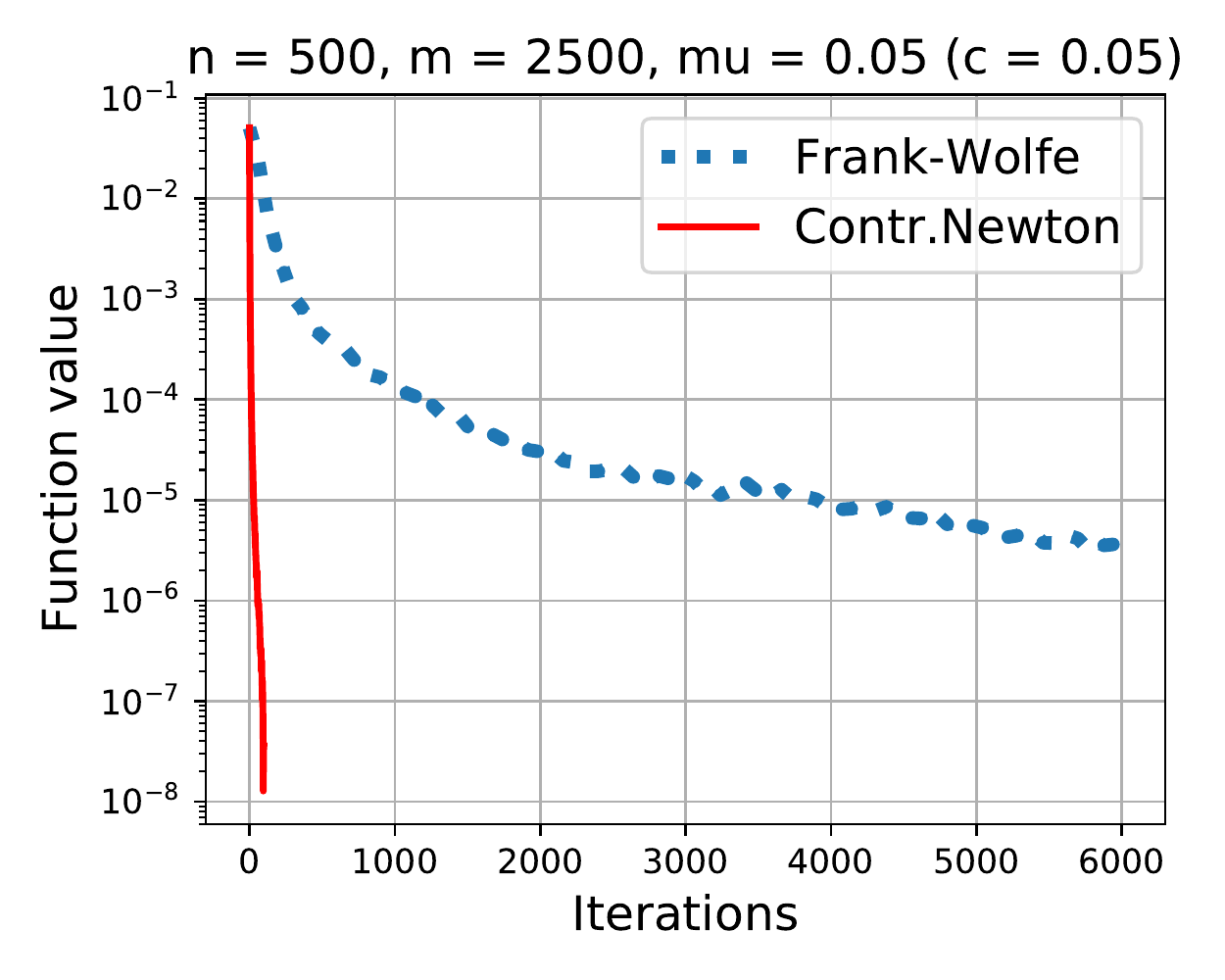}
	\includegraphics[width=0.21\textwidth ]{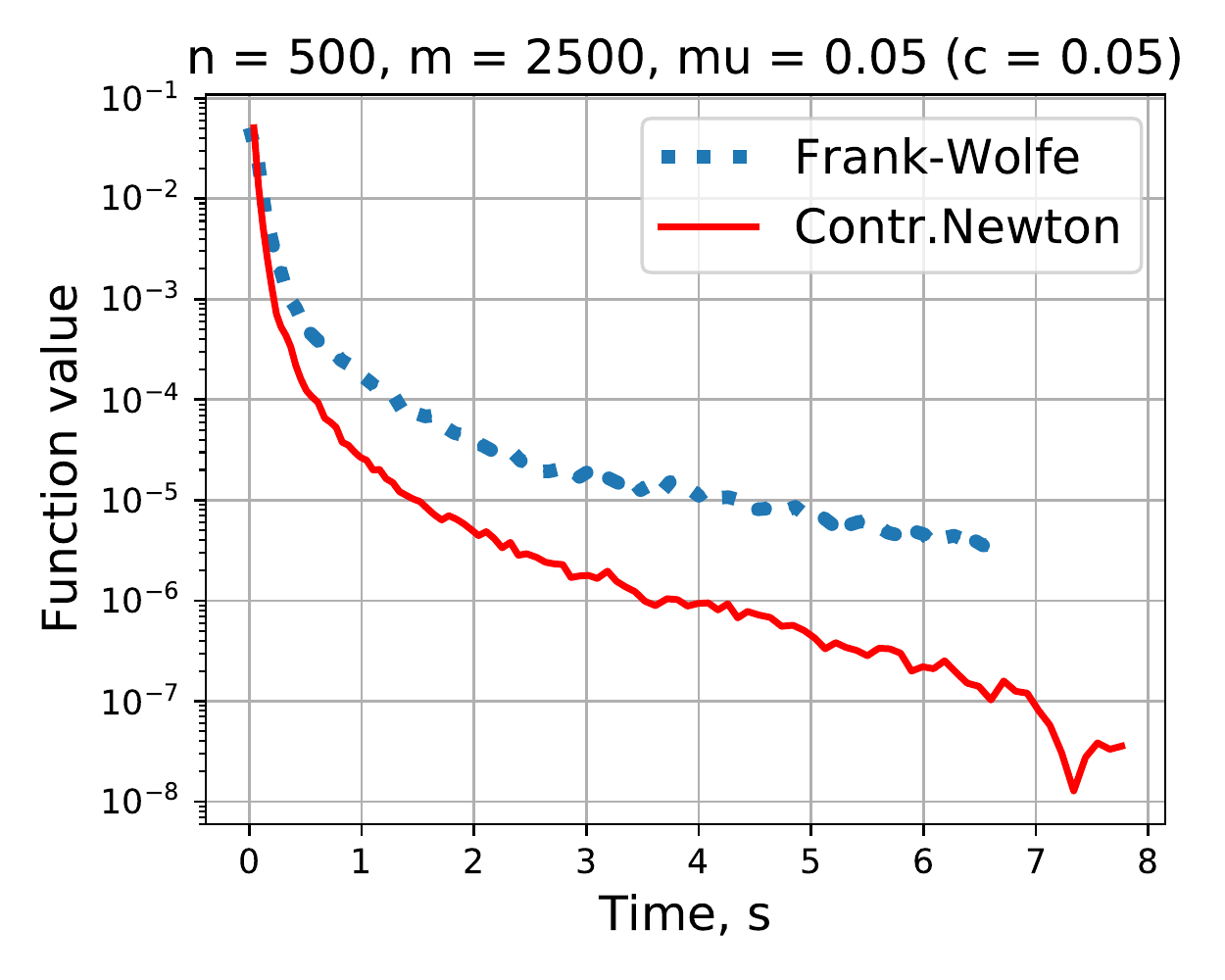}
	\caption{$n = 500, \; m = 2500$.}
	\label{Fig3}
\end{figure}

Let us consider the problem of minimizing the log-sum-exp function (SoftMax)
$$
\ba{rcl}
f_{\mu}(x) & = & 
\mu \ln \biggl( \sum\limits_{i = 1}^m \exp\Bigl(  \frac{\la a_i, x \ra - b_i}{\mu} \Bigr) \biggr),
\qquad x \in \R^n,
\ea
$$
over the standard simplex $\S_n$~\eqref{SimplexDef}.
Coefficients $\{ a_i  \}_{i = 1}^m$ and $b$ are generated randomly from the uniform
distribution on $[-1, 1]$.
We compare the performance of Inexact Contracting Newton Method~\eqref{met-CPHF} 
with that one of the classical Frank-Wolfe algorithm,
for different values of the parameters.
The results are shown on Figures~\ref{Fig1}-\ref{Fig3}.

We see, that the new method
works significantly better in terms of the outer iterations (oracle calls).
This confirms our theory.
At the same time, for many values of the parameters, it shows better performance
in terms of total computational time as well\footnote{
	Clock time was evaluated using the machine with Intel Core i5 CPU, 1.6GHz; 8 GB RAM. 
	The methods were implemented in Python.}.

\newpage

\section{Discussion}
\label{sc-Discussion}

In this paper, we present a new general framework of Contracting-Point methods,
which can be used for developing affine-invariant optimization algorithms 
of different order. 
For the methods of order $p \geq 1$,
we prove the following global convergence rate:
$$
\ba{rcl}
F(x_k) - F^{*} & \leq & {\cal O}\bigl( 1 / k^{p} \bigr), \quad k \geq 1.
\ea
$$
This is the same rate, as that of the basic 
high-order Proximal-Point scheme~\cite{nesterov2020inexact}.
However, the methods from our paper are free from using the norms
or any other characteristic parameters of the problem. 
This nice property makes Contracting-Point methods favourable 
for solving optimization problems over the sets with a non-Euclidean geometry
(e.g. over the simplex or over a general convex polytope).

At the same time, it is known that in Euclidean case, the prox-type methods
can be accelerated, achieving 
${\cal O}(1 / k^{p + 1})$ global rate of convergence~\cite{baes2009estimate,nesterov2019implementable,doikov2019contracting,nesterov2020inexact}.
Using additional one-dimensional search at each iteration, this rate can be improved up to 
${\cal O}(1 / k^{\frac{3p + 1}{2}})$ (see~\cite{gasnikov2019near,nesterov2020inexact}).
The latter rate is shown to be optimal~\cite{arjevani2019oracle,nesterov2018lectures}.
To the best of our knowledge, the lower bounds for high-order methods in general non-Euclidean case
remain unknown. However, the worst-case oracle complexity of the classical Frank-Wolfe algorithm 
(the case $p = 1$ in our framework) is proven to be near-optimal
for smooth minimization over $\| \cdot \|_{\infty}$-balls~\cite{guzman2015lower}.
 
Another open question is a possibility of efficient implementation of our methods for the case $p \geq 3$.
In view of absence of explicit regularizer (contrary to the prox-type methods), the subproblem in~\eqref{TensorSub}
can be nonconvex. Hence, it seems hard to find its global minimizer.
We hope that for some problem classes, it is still feasible to satisfy the inexact stationarity condition~\eqref{TensorStationary} by reasonable amount of computations.
We keep this question for further investigation.


\section*{Appendix}

In this section, we state some 
simple facts about maximization of multilinear symmetric forms.
For a fixed $p \geq 1$,
let us consider $(p + 1)$-linear symmetric form
$A$. For a set of vectors $h_1, \dots, h_{p + 1} \in \E$, we have
$$
\ba{rcl}
A[h_1, \dots, h_{p + 1}] & \in & \R.
\ea
$$
For two vectors $u, v \in \E$
and integers $i, j \geq 0$ such that $i + j = p + 1$,
we use the following shorter notation:
$$
\ba{rcl}
A[u]^{i}[v]^{j}
& \Def &
A[\underbrace{u, \dots, u}_{i \; \text{times}}, \underbrace{v, \dots, v}_{j \; \text{times}}].
\ea
$$
Let us fix arbitrary compact convex set $S \subset \E$.
We are interested to bound  the variation of $A$
over two vectors from $S$, by that over the only one vector:
\beq \label{VarBound}
\ba{rcl}
\sup\limits_{u, v \in S} | A[u]^{p}[v] | & \leq &
{\mathcal C}_p  \sup\limits_{h \in S} | A[h]^{p + 1}|,
\ea
\eeq
for some constant ${\mathcal C}_p$.
Note, that if $S$ is a ball in the Euclidean norm,
then ${\mathcal C}_p = 1$, and
the values of both supremums are equal
(see Appendix~1 in~\cite{nesterov1994interior}, 
and Section 2.3 in~\cite{nemirovski2004interior}).
In what follows, our aim is to estimate the value of ${\mathcal C}_p$ 
for arbitrary $S$. Namely, we establish the following bound.

\BP \label{propMultilinear}
For any compact convex set $S$, \eqref{VarBound} holds
with
\beq \label{CpBound}
\ba{rcl}
{\mathcal C}_p & = & \frac{ (p + 1)^{p + 1} + p^{p + 1} + 1 }{(p + 1)!}
\;\; \leq \;\; 
\frac{ 2(p + 1)^p }{p!}.
\ea
\eeq
\EP
\proof
For a pair of integers $n, k \geq 0$, let us denote by 
$\binom{n}{k}$ the binomial coefficients, given by the formula
$$
\ba{rcl}
\binom{n}{k} & \Def & \frac{n  (n - 1) \cdots (n - k + 1)}{k!},
\ea
$$ 
and by ${n \brace k}$ we denote the Stirling numbers of the second kind.
By definition, ${n \brace k}$ is equal to the number of ways
to partition a set of $n$ objects into $k$ nonempty subsets. 
The following important identity holds 
(see, for example, \cite{graham1989concrete}):
\beq \label{PowerStirling}
\ba{rcl}
k^n & = & k! \sum\limits_{r = 1}^{k} \frac{ {n \brace r} }{(k - r)!}, \qquad n \geq 1.
\ea
\eeq
Note also, that ${n \brace n} = 1$, and ${n \brace k} = 0$, for $k > n$. 

Now, let us fix arbitrary vectors $u, v \in S$,
and consider the set of their convex combinations $h_i = \alpha_i u + (1 - \alpha_i) v \in S$
for some $\alpha_i \in (0, 1)$, $1 \leq i \leq p$. The binomial theorem yields the system of equations
\beq \label{PolySystem}
\ba{rcl}
A[h_i]^{p + 1} & = & 
\sum\limits_{j = 0}^{p + 1} 
\binom{p + 1}{j} \alpha_i^{j} (1 - \alpha_i)^{p + 1 - j} A[u]^j [v]^{p + 1 - j},
\qquad 1 \leq i \leq p.
\ea
\eeq
For the choice $\alpha_i = \frac{i}{i + 1}$, we have $1 - \alpha_i = \frac{1}{i + 1}$,
and
$$
\ba{rcl}
\alpha_i^{j} (1 - \alpha_i)^{p + 1 - j} & = & \frac{i^j}{(i + 1)^{p + 1}}.
\ea
$$
Therefore, introducing a vector $x \in \R^p$,
$$
\ba{rcl}
x^{(j)} & \equiv & \binom{p + 1}{j}A[u]^j[v]^{p + 1 - j}, \qquad 1 \leq j \leq p,
\ea
$$
from~\eqref{PolySystem} we obtain the linear system $\boxed{Bx = c}$ with matrix 
\beq \label{PolySystemMat}
\ba{rcl}
B^{(i, j)} & \equiv & i^j, \qquad 1 \leq i, j \leq p,
\ea
\eeq
and the right hand side vector
\beq \label{PolySystemC}
\ba{rcl}
c^{(i)} & \equiv & 
(i + 1)^{p + 1} \bigl(  
A[h_i]^{p + 1} - (1 - \alpha_i)^{p + 1} A[v]^{p + 1} - \alpha_i^{p + 1} A[u]^{p + 1} \bigr),
\qquad 1 \leq i \leq p.
\ea
\eeq
The matrix given by~\eqref{PolySystemMat} looks as follows
$$
\ba{rcl}
B & = & \begin{pmatrix}
&1 &1 &1 &\dots &1 \\
&2 &2^2 &2^3 &\dots &2^p \\
&3 &3^2 &3^3 &\dots &3^p \\
&\vdots &\vdots &\vdots &\ddots &\vdots \\
&p &p^2 &p^3 &\dots &p^p
\end{pmatrix}.
\ea
$$
This structure is similar to that one of the Vandermonde matrix.
By the Gaussian elimination process we can build a sequence of matrices
$$
\ba{rcl}
B & = & B_1  \;\; \mapsto \;\; B_2 \;\; \mapsto \;\; \dots \;\; \mapsto \;\; B_p,
\ea
$$
such that $B_p$ is upper triangular, and the corresponding sequence 
of the right hand side vectors
$$
\ba{rcl}
c & = & c_1 \;\; \mapsto \;\; c_2 \;\; \mapsto \;\; \dots \;\; \mapsto \;\; c_p,
\ea
$$
having the same solution $x$ as the initial system:
$$
\ba{rcl}
B_t x = c_t, \qquad 1 \leq t \leq p.
\ea
$$
Then, the last component of the solution can be easily found:
\beq \label{AuvEqual}
\ba{rcl}
(p + 1) A[u]^p[v] & = &  x^{(p)} \;\; = \;\; \frac{c_p^{(p)}}{B_p^{(p, p)}},
\ea
\eeq
from which we may obtain the required bound for the left hand side of~\eqref{VarBound}.
Thus, we are interested to investigate the elements of $B_p$ and $c_p$.

Let us prove by induction, that for every $1 \leq t \leq p$, it holds
\beq \label{PolySystemGauss}
\ba{rcl}
B_t^{(i, j)} & = & \begin{cases}
	i! {j \brace i} & \text{if} \; i \leq t; \\
	i! \sum\limits_{r = t }^i \frac{ {j \brace r}  }{(i - r)!} & \text{otherwise}.
\end{cases}
\ea
\eeq
For $t = 1$,~\eqref{PolySystemGauss} follows from~\eqref{PowerStirling}, 
and this is the base of the induction.
At step $t$ of the Gaussian elimination, we have the matrix $B_t$.
First, we freeze its $t$-th row for all the following matrices:
$$
\ba{rcl}
B_t^{(t, j)} & = & B_{t + 1}^{(t, j)} 
\;\; = \;\; B_{t + 2}^{(t, j)}
\;\; = \;\; \dots \;\; = \;\;
B_{p}^{(t, j)} \;\; = \;\; t! { j \brace t }, \qquad 1 \leq j \leq p.
\ea
$$
Then, we subtract this row from all the rows located below,
scaled by an appropriate factor,
for $t < i \leq p$:
$$
\ba{rcl}
B_{t + 1}^{(i, j)} & = & B_t^{(i, j)} -  \frac{B_t^{(i, t)}}{B_t^{(t, t)}} \cdot B_t^{(t, j)}, \qquad 1 \leq j \leq p.
\ea
$$
Note, that $B_t^{(t, t)} = t!$ and $B_t^{(i, t)} = \frac{i!}{(i - t)!}$. Therefore, we obtain
$$
\ba{rcl}
B_{t + 1}^{(i, j)}
& = & 
i! \sum\limits_{r = t}^{i} \frac{ {j \brace r} }{(i - r)!}
\;\; - \;\;
\frac{i! {j \brace t}}{(i - t)!}
\;\; = \;\;
i! \sum\limits_{r = t + 1}^{i} \frac{ {j \brace r} }{(i - r)!},
\qquad 1 \leq j \leq p,
\ea
$$
and this is~\eqref{PolySystemGauss} for the next step. 
Hence~\eqref{PolySystemGauss} is established by induction for all $1 \leq t \leq p$.

Similarly, we have the update rules for the right hand sides:
$$
\ba{rcl}
c_t^{(t)} & = & c_{t + 1}^{(t)}
\;\; = \;\; \dots 
\;\; = \;\;
c_p^{(t)},
\ea
$$
and for $t < i \leq p$:
$$
\ba{rcl}
c_{t + 1}^{(i)} & = & c_t^{(i)} - \frac{B_t^{(i, t)}}{B_t^{t, t}} c_t^{(t)}
\;\; = \;\; c_t^{(i)} - \binom{i}{t} c_t^{(t)}  \\
\\
& = & c_{t - 1}^{(i)} - \binom{i}{t - 1} c_{t - 1}^{(t - 1)} - \binom{i}{t} c_t^{(t)} \;\; = \;\; \dots \\
\\
& = & c_1^{(i)} - \sum\limits_{r = 1}^{t} \binom{i}{r} c_r^{(r)}.
\ea
$$
Therefore, for the resulting vector $c_p$ we have a recurrence:
\beq \label{CRecurr}
\ba{rcl}
c_p^{(i)} = c_1^{(i)} - \sum\limits_{r = 1}^{i - 1} \binom{i}{r} c_p^{(r)},
\qquad 1 \leq i \leq p.
\ea
\eeq
From~\eqref{CRecurr} we obtain an explicit expression for $c_p$ using only the initial values:
\beq \label{CExplicit}
\ba{rcl}
c_p^{(i)} & = & \sum\limits_{j = 1}^i (-1)^{i - j} \binom{i}{j} c^{(j)}_1,
\qquad 1 \leq i \le p.
\ea
\eeq
Indeed, for $i = 1$~\eqref{CExplicit} follows directly from~\eqref{CRecurr}.
Assume by induction that~\eqref{CExplicit} holds for all $1 \leq i \leq n$,
for some $n$.
Then, for the next index we have
$$
\ba{rcl}
c_p^{(n + 1)} & \refEQ{CRecurr} & 
c_1^{(n + 1)} - \sum\limits_{r = 1}^n \binom{n + 1}{r} c_p^{(r)} \\
\\
& \refEQ{CExplicit} &
c_1^{(n + 1)}
- \sum\limits_{r = 1}^n \sum\limits_{j = 1}^r (-1)^{r - j} \binom{n + 1}{r} \binom{r}{j} c_1^{(j)} \\
\\
& = &
c_1^{(n + 1)} - \sum\limits_{j = 1}^n  
\biggl(\sum\limits_{r = j}^{n} (-1)^{r - j}\binom{n + 1}{r}\binom{r}{j} \biggr) c_1^{(j)} \\
\\
& = &
c_1^{(n + 1)} + \sum\limits_{j = 1}^n  (-1)^{n + 1 - j}\binom{n + 1}{j} c_1^{(j)},
\ea 
$$
where the last equation follows from simple observations:
$$
\ba{rcl}
\sum\limits_{r = j}^{n} (-1)^{r - j}\binom{n + 1}{r}\binom{r}{j}
& = & 
\sum\limits_{r = j}^{n} (-1)^{r - j} \frac{(n + 1)! \, r!}{r! \, (n + 1 - r)! \, j! \, (r - j)!} \\
\\
& = &
\binom{n + 1}{j} \sum\limits_{r = j}^{n} (-1)^{r - j} \binom{n + 1 - j}{r - j} \\
\\
& = & 
\binom{n + 1}{j} \sum\limits_{l = 0}^{n - j} (-1)^{l} \binom{n + 1 - j}{l} \\
\\
& = &
\binom{n + 1}{j} \Bigl( (1 - 1)^{n + 1 - j} - (-1)^{n + 1 - j}  \Bigr) 
\;\; = \;\; (-1)^{n - j} \binom{n + 1}{j}.
\ea
$$
Hence~\eqref{CExplicit} is established by induction for all $1 \leq i \le p$.

Let us denote by ${\mathcal{V}}$ the supremum from the right hand side of~\eqref{VarBound}:
$$
\ba{rcl}
\mathcal{V} & \Def & \sup\limits_{h \in S} |A[h]^{p + 1}|.
\ea
$$
Then, in view of~\eqref{PolySystemC}, we have
\beq \label{C1Bound}
\ba{rcl}
| c_1^{(j)} | 
& = &
| c^{(j)} | 
\;\; \leq \;\;
((p + 1)^{p + 1} + p^{p + 1} + 1) \mathcal{V}, \qquad 1 \leq j \leq p,
\ea
\eeq
and, consequently
$$
\ba{rcl}
| A[u]^p[v] |
& \overset{\eqref{AuvEqual},\eqref{PolySystemGauss}}{=} & 
\frac{|c_p^{(p)}|}{(p + 1)!}
\;\; \overset{\eqref{CExplicit}}{=} \;\;
\frac{1}{(p + 1)!} 
\Big|  \sum\limits_{j = 1}^p (-1)^{p - j} \binom{p}{j} c_1^{(j)}  \Big| \\
\\
& \overset{\eqref{C1Bound}}{\leq} &
\frac{( (p + 1)^{p + 1} + p^{p + 1} + 1  ) {\mathcal{V}}}{(p + 1)!}
\Big| \sum\limits_{j = 1}^p (-1)^{p - j} \binom{p}{j}  \Big|
\;\; = \;\;
\mathcal{C}_p \mathcal{V}.
\ea
$$
Since $u, v \in S$ are arbitrary vectors, we have~\eqref{CpBound} established.
\qed

Let us consider the most important cases, when $p = 1$ and $p = 2$.
\BC
For any symmetric bilinear form $A: \E \times \E \to \R$ and
any compact convex set $S \subset \E$, it holds
\beq \label{ColBilinear}
\ba{rcl}
\sup\limits_{u, v \in S} |A[u, v]| & \leq & 
3 \sup\limits_{h \in S} | A[h, h] |.
\ea
\eeq
\EC

\BC
For any symmetric trilinear form $A: \E \times \E \times \E \to \R$ and
any compact convex set $S \subset \E$, it holds
\beq \label{ColTrilinear}
\ba{rcl}
\sup\limits_{u, v \in S} |A[u, u, v]| & \leq & 
6 \sup\limits_{h \in S} | A[h, h, h] |.
\ea
\eeq
\EC

It appears that the bound in~\eqref{ColBilinear} is tight.
\BE
Consider the following symmetric bilinear form on two-dimensional space $\E = \R^2$:
$$
\ba{rcl}
A[u, v] & = & u^{(1)} v^{(1)} - 2 u^{(2)} v^{(2)}, \qquad u, v \in \R^{2},
\ea
$$
and let
$$
\ba{rcl}
S & = & \bigl\{  x \in \R^2 \; : \; x^{(1)} = 1, \;   x^{(2)} \in [-1, 1]  \bigr\}.
\ea
$$
Then,
$$
\ba{rcl}
\sup\limits_{u, v \in S} |A[u, v]| & = & 
\sup\limits_{\alpha, \beta \in [-1, 1] } |1 - 2 \alpha \beta | 
\;\; = \;\; 3.
\ea
$$
However,
$$
\ba{rcl}
\sup\limits_{h \in S} |A[h,  h]|
& = &
\sup\limits_{\alpha \in [0, 1]} |1 - 2\alpha | 
\;\; = \;\; 1.
\ea
$$
\EE 

If it happens that our bilinear form is positive semidefinite 
(e.g. it is determined by the Hessian of a convex function),
the constant in~\eqref{ColBilinear} can be improved to be $1$,
so the both supremums are equal.

\BP Let symmetric bilinear form $A: \E \times \E \to \R$ be 
positive semidefinite:
$$
\ba{rcl}
A[h, h] & \geq & 0, \qquad h \in \E.
\ea
$$
Then, for \underline{any} set $S \subset \E$, it holds
$$
\ba{rcl}
\sup\limits_{u, v \in S} |A[u, v]|
& = & \sup\limits_{h \in S} A[h, h].
\ea
$$
\EP
\proof
Indeed, by the Eigenvalue Decomposition, for some $r \geq 0$,
there exists a set of linear forms $a_1, \dots, a_r \in \E^{*}$
and positive numbers $\lambda_1, \dots, \lambda_r > 0$ such that
$$
\ba{rcl}
A[u, v] & = & \sum\limits_{i = 1}^r \lambda_i \la a_i, u \ra \la a_i, v \ra, \qquad u, v \in S.
\ea
$$
Therefore, using Cauchy-Bunyakovsky-Schwarz inequality, we get
$$
\ba{rcl}
|A[u, v]| & \leq & 
\Bigl(  \sum\limits_{i = 1}^r \lambda_i \la a_i, u\ra^2  \Bigr)^{1/2} 
\Bigl(  \sum\limits_{i = 1}^r \lambda_i \la a_i, v \ra^2 \Bigr)^{1/2} \\
\\
&=& \bigl(  A[u, u]  \bigr)^{1/2} \bigl(  A[v, v]  \bigr)^{1/2} \\
\\
& \leq & 
\sup\limits_{h \in S} A[h, h].
\ea
$$

\qed

\end{document}